\documentclass[11pt]{article}
\pdfoutput=1
\usepackage[a4paper,margin=2.5cm]{geometry}
\usepackage{amsmath,amssymb}
\usepackage{dsfont}
\usepackage{accents}
\usepackage{wrapfig}

\usepackage{graphicx}
\usepackage{url}
\newcommand{\R}{\ensuremath{\mathbb{R}}}

\newcommand{\Hyp}{\ensuremath{\mathbb{H}}}

\newcommand{\M}{\ensuremath{{\cal M}}}

\newcommand{\U}{\ensuremath{{\cal U}}}
\newcommand{\Sph}{\ensuremath{{\cal S}}}

\newcommand{\MM}{\ensuremath{{\mathfrak M}}}
\newcommand{\NN}{\ensuremath{{\mathfrak N}}}
\newcommand{\XX}{\ensuremath{{\mathfrak X}}} %\widehat {\mathfrak M}
\newcommand{\contract}{\ensuremath{\: \overset{\scriptscriptstyle{\bullet}}{{\scriptscriptstyle{\bullet}}} \: }}
\newcommand{\contractvar}{\ensuremath{{\scriptscriptstyle{\bullet}}}}
\newcommand{\contractcirc}{\ensuremath{\: \overset{\scriptscriptstyle{\circ}}{{\scriptscriptstyle{\circ}}} \: }}
\newcommand{\contractvarcirc}{\ensuremath{{\scriptscriptstyle{\circ}}}}

\newcommand{\Cut}{\ensuremath{\text{Cut}}}
\newcommand{\Inj}{\ensuremath{\text{Inj}}}
\newcommand{\inj}{\ensuremath{\text{inj}}}

\newcommand{\dist}{\ensuremath{\:\mbox{\rm dist}}}
\newcommand{\scal}[2]{\ensuremath{ \left<
    \:#1\:\left|\:#2\right.\right> }}
\newcommand{\inv}{^{\text{\tiny (-1)}}}
\newcommand{\expect}[1]{\ensuremath{ \mathbf{E}\left[\: #1\:\right]}}

\newcommand{\Id}{\ensuremath{\:\mathrm{Id}}}

\newcommand{\lcp}[1]{\ensuremath{\overrightarrow{#1}}}

\newtheorem{theorem}{Theorem}
\newtheorem{definition}{Definition}

\begin{document}
\pagenumbering{arabic}
\pagestyle{headings} 

\title{Curvature effects on the empirical mean in Riemannian and affine Manifolds:
a non-asymptotic  high concentration expansion in the small-sample regime}

\author{Xavier Pennec\footnote{Universti\'e C\^ote d'Azur and Inria, Sophia Antipolis, France. Email: xavier.pennec@inria.fr}}

\maketitle

\begin{abstract}
The asymptotic concentration of the Fr\'echet mean of IID random variables on a Riemannian manifold was established with a central limit theorem by Bhattacharya \& Patrangenaru (BP-CLT) \cite{bhattacharya_large_2005}. This asymptotic result shows that the Fr\'echet mean behaves almost as the usual Euclidean case for sufficiently concentrated distributions. However, the asymptotic covariance matrix of the empirical mean is modified by the expected Hessian of the squared distance.  This Hessian matrix was explicitly computed in \cite{bhattacharya_statistics_2008} for constant curvature spaces in order to relate it to the sectional curvature. Although explicit, the formula remains quite difficult to interpret, and the intuitive effect of the curvature on the asymptotic convergence remains unclear. Moreover, we are most often interested  in the mean of a finite sample of small size in practice. In this work, we aim at understanding the effect of the manifold curvature in this small sample regime. Last but not least, one would like computable and interpretable approximations that can be extended from the empirical Fr\'echet mean in Riemannian manifolds to the empirical exponential barycenters in affine connection manifolds. 

For distributions that are highly concentrated around their mean, and for any finite number of samples,  we establish explicit Taylor expansions on the first and second moment  of the empirical mean 
thanks to a new Taylor expansion of the Riemannian log-map in affine connection spaces. 
This shows that the empirical mean has a bias in $1/n$ proportional to the gradient of the curvature tensor contracted twice with the covariance matrix, and a modulation of the convergence rate of the covariance matrix proportional to the covariance-curvature tensor. 
We show that our non-asymptotic high concentration expansion is consistent with the asymptotic expansion of the BP-CLT.  Experiments on constant curvature spaces demonstrate that both expansions are very accurate in their domain of validity. Moreover, the modulation of the convergence rate of the empirical mean's covariance matrix is explicitly encoded using a scalar multiplicative factor that gives an intuitive vision of the impact of the curvature: the variance of the empirical mean decreases faster than in the Euclidean case in negatively curved space forms, with an infinite speed for an infinite negative curvature. This suggests potential links with the stickiness of the Fr\'echet mean described in stratified spaces.
On the contrary, the variance of the empirical mean decreases  more slowly  than in the Euclidean case in positive curvature space forms, with divergence when we approach the limits of the Karcher \& Kendall concentration conditions with a uniform distribution on the equator of the sphere, for which the Fr\'echet mean is not a single point any more. 
\end{abstract}

%%%%%%%%%%%%%%%%%%%%%%%%%%%%%%%%%%%%%%%%%%%%%%%%%%%%%%%%%%%%%%%%%%%%%%%%%
\clearpage

\section{Introduction and overview of the results}
The asymptotic concentration of the Fr\'echet mean of IID random variables on a Riemannian manifold was established in a central limit theorem by Bhattacharya \& Patrangenaru (BP-CLT) \cite{bhattacharya_large_2005}. This asymptotic result showed that the Fr\'echet mean behaves almost as the usual Euclidean case  for sufficiently concentrated distributions once we have taken into account the curvature effects in the Hessian of the variance.  These results were later generalized in \cite{kendall_limit_2011} to establish weak laws of large numbers and central limit theorems of Lindeberg type for empirical Fr\'echet means of independent but non-identically distributed random variables. The essential quantity modifying the asymptotic covariance matrix of the empirical mean is the Hessian of the squared variance (the expectation of the Hessian of the squared distance). 
This Hessian matrix was explicitly computed in \cite{bhattacharya_statistics_2008} for constant curvature spaces in order to relate the asymptotic dispersion of the Fr\'echet mean  to the sectional curvature. Although explicit, the formula remains quite difficult to interpret, and the intuitive effect of the curvature on the asymptotic convergence remains cryptic. Moreover, we most often have in practice a finite sample of small size and it is interesting to understand if the manifold curvature has addition effects on the estimation of the mean in this small sample regime. Last but not least, we would like to obtain results that can be extended from the Fr\'echet mean in Riemannian manifolds to exponential barycenters in affine connection manifolds.

We investigate in this work  the moments of the distribution of the empirical mean of a fixed number of samples in Riemannian and affine connection manifolds. We derive interpretable approximations for sufficiently concentrated distributions showing how the moments of the empirical mean deviate from the usual Euclidean case in the non-asymptotic case. We also study the speed of convergence towards their asymptotic value.
In order to obtain a high concentration expansion of the minimum of the expected Riemannian squared distance, we may base our computations on a Taylor expansion of the Riemannian metric in a normal coordinate system at a fixed point.
Such an extension is well known for the first orders, but going to orders higher than 4 is computationally much more involved and requires computational algebraic methods such as the ones developed in \cite{brewin_riemann_2009}. This leads to quite involved formulas with contractions over many covariant and contravariant indices\footnote{The results presented in this paper were originally developed in 2015 using Brewin formulas of this type. However, the difficulty to make sense of pages of computations involving summations over many indices led us to delay the publication until we found a simpler and more intuitive coordinate free formulation.}.  
In this paper, we base our Taylor expansions on a coordinate free expansion of the composition of two exponential maps developed by Gavrilov \cite{gavrilov_algebraic_2006,gavrilov_double_2007}. This formula only involves the curvature and torsion tensors and their covariant derivatives at the development point. Moreover, the proposed extension relies only on the connection and not on the Riemannian metric, which makes it also suitable for the more general case of affine connection spaces with distributions sufficiently concentrated around their exponential barycenter.

\paragraph{Assumptions and objectives}

Let $X_n = \sum_{i=1}^n \delta_{x_i}$ be an IID $n$-sample drawn from a probability distribution $\mu$ on a Riemannian manifold $\M$.  We want to investigate how the moments of the empirical Fr\'echet mean $\bar x_n$ of the empirical distribution $X_n$ compares to the ones of the population Fr\'echet mean $\bar x$ of $\mu$.  This formulation obviously assumes that the Fr\'echet mean is unique. Thus, we assume that the distribution $\mu$ belongs to a sufficiently small strongly convex neighborhood $\U$ of diameter $\varepsilon$ in $\M$. In Riemannian manifolds, this can be quantified by the Kendall and Karcher concentration (KKC) condition. 
This allows us to lift the distribution and all the computations to the convex subset $V_x = \log_x(\U) \subset T_x\M$ of the tangent space at any point $x \in \U$. 
In a convex affine connection space, the notion of mean that makes sense in convex subsets is based on exponential barycenters. One can show uniqueness in Arnaudon \& Li convexity (ALC) conditions for convex affine manifold with p-convex geometry. This notion involves an auxiliary metric with which we can compute the diameter $\varepsilon = \text{diam}(\U)$ of the convex set $\U$. 

Note that any   distribution $\mu$ with a support in $\U$, including the empirical distribution $X_n = \frac{1}{n}\sum_{i=1}^n \delta_{x_i}$, has a unique mean.  Moreover, the $k$-th order moments $\MM_k\mu(x) =  \int_{\M}\log_x(z)^{\otimes k} \mu(dz)$ is of order $k$ in $\varepsilon = \text{diam}(\U)$ since $\| \log_x(z)\|_x \leq \varepsilon$ for any $x$ and $z \in V_x$. 
For our expansions, we also need to control the order of curvature terms that will appear. The Riemannian curvature tensor $R(u,v)w$ being multilinear in its variables, it is a term of order 3 in $\varepsilon$ at each point $x\in \M$ for vectors $u,v,w \in V_x = \log_x( \U ) \subset B(x,\varepsilon) \subset T_x\M$. Likewise, the covariant derivative $\nabla_t R(u,v)w$ is of order 4 for vectors in the same subset.

A distribution $\mu$ satisfying the KKC or the ALC conditions has a unique mean in $\U$ implicitly defined with the exponential barycenter formulation: 
$\MM_1\mu(z) =  \int_{\M}\log_z(y) \: \mu(dy)=0$.
The goal is to find a high order approximation of the solution of the barycentric equation $\MM_1\mu(z) =0$ in the neighborhood of a point $x \in \U$. For that purpose, we parametrize the points of $\U$ by $x_v = \exp_x(v)$ with $v \in V_x = \log_x(\U) \subset T_x\M$ and we look for a Taylor expansion of $\MM_1\mu(x_v)=0$ with respect to $v$. However, the tangent mean $\MM_1\mu(x_v) = \int_\U \log_{x_v}(y)\: \mu(dy)$ is a vector field on $\U$, i.e. a mapping from $\U \subset \M$ to $T_{x_v}\M$. In order to have a unique image space, we may take a chart, which would imply working in a specific coordinate system, or more interestingly parallel transport each tangent vector from $x_v$ to $x$: we obtain this way the recentered tangent mean map (Definition \ref{def:RecenterdMeanMap}): 
\[
\NN_x^{\mu}(v) = \Pi_{\exp_x(v)}^x \MM_1\mu(\exp_x(v)) \in  T_x\M.
\] 
The recentered tangent mean map is a mapping of vector spaces from $V_x\subset T_x\M$ to $T_x\M$ whose zeros parametrize the exponential barycenters of $\mu$. In order to localize them, we aim at computing a series  expansion of the recentered tangent mean map with respect to $v$.

\paragraph{Taylor expansions in affine connection spaces}
For that purpose, we develop methods for Taylor expansions in manifolds in Section \ref{sec:Taylor}. Based on the coordinate free expansion of the composition of two exponential maps of \cite{gavrilov_algebraic_2006,gavrilov_double_2007} (the double exponential, theorem \ref{thm:DoubleExpExpansion}), we derive  a series expansion up to order 5 of the logarithm of a fixed point $x_w=\exp_x(w)$ at a point $x_v=\exp_x(v)$, parallel transported back to $x$ (Theorem \ref{thm:NLogExpansion}): 
\[
\begin{split}
l_x(v, w) = &  \Pi_{x_v}^x \log_{x_v}(\exp_x(w)) = w-v + \frac{1}{6}R(w ,v)(v-2 w)   \\
&+ \frac{1}{24}(\nabla_v R)(w,v) (2v-3w) 
 + \frac{1}{24}(\nabla_{w} R)(w,v) (v-2w) +O(5).
\end{split}
\]
This neighboring log expansion is non-metric and valid for general affine connection manifolds. 
In the Riemannian case, we can take the square norm of that vector to get an expansion of the square Riemannian distance between two points that are away from the base point $x$. This expansion is the one needed to write an expansion of the variance around the point $x$. 

However, since the minimum of the variance is in particular a critical point and thus a zero of the recentered tangent mean map $\NN_x^{\mu}(v)$, it is more convenient to compute directly the polynomial expansion of this vector space mapping. This is the focus of Section \ref{sec:FrechetMeanExpansion}. Thanks to the previous work, it is relatively straightforward to show that
 \[
%\label{eq:RecenteredMeanExpansion}
\begin{split}
\NN_x^{\mu}(v) = & \: \MM_1  - v
+ \frac{1}{6} R(\MM_1,v)v  -\frac{1}{3}R(\contractvar,v)\contractvar \contract \MM_2 
 +\frac{1}{12} (\nabla_v R)( \MM_1,v)v \\
& + \frac{1}{24} (\nabla_\contractvar R) (\contractvar,v)v \contract \MM_2 
 -\frac{1}{8} (\nabla_v R)(\contractvar,v)\contractvar \contract \MM_2 -\frac{1}{12} (\nabla_\contractvar R)(\contractvar, v)\contractvar \contract \MM_3 + O(\varepsilon^5)
\end{split}
\]
where the notation $R(\contractvar, v)\contractvar \contract \MM_2$ denotes the contraction of the tensor moment with the curvature tensor $R$ along the axes specified by the bullets.
 
Solving for the value of $v = \log_x(\bar x)$ that zeros out this expression leads to the polynomial expansion of the field $\log_x(\bar x)$ pointing from the points $x$ to the mean $\bar x$ (theorem \ref{thm:FrechetMeanField}):
 \[
% \label{LogFrechet}
\begin{split}
\log_x(\bar x) = &
\MM_1  - \frac{1}{3} R(\contractvar, \MM_1)\contractvar \contract \MM_2 
- \frac{1}{24} \nabla_\contractvar R(\contractvar,\MM_1)\MM_1 \contract \MM_2 \\ 
& -\frac{1}{8} \nabla_{\MM_1} R(\contractvar, \MM_1) \contractvar \contract \MM_2
-\frac{1}{12} \nabla_{\contractvar} R(\contractvar, \MM_1) \contractvar \contract \MM_3 + O(\epsilon^5).
\end{split}
\] 

\paragraph{Non asymptotic high concentration expansion of the moments of the empirical mean}
Equipped with this expansion, we analyze in Section \ref{sec:IIDSample} the expected first moment $\expect{\log_x(\bar x_n) }$ and the expected second moment  $\expect{\log_x(\bar x_n) \otimes \log_x(\bar x_n)}$ of the empirical Fr\'echet mean $\bar x_n$ of $n$ IID samples at the population mean $\bar x$. Using the previous expansion gives us a formula involving the tensor product of empirical moments of a sample. Taking the expectation of empirical moments is simple: it gives the moment of the underlying distribution. Taking  the expectation of tensor products of empirical moments is more complex since it generates some tensor products of moments of mixed orders. 
We establish in Theorem \ref{thm:MomentsEmpiricalMean} that the expected log of the empirical mean at the population mean is:
\[
%\text{Bias}(\bar x_n) = 
\expect{ \log_{\bar x}(\bar x_n) } =  
  \textstyle  \frac{1}{6 n}\left(1-\frac{1}{n}\right)  \MM_2 \contractcirc \nabla_{\contractvar} R(\contractvar, \contractvarcirc) \contractvarcirc \contract \MM_2  
 + O\left(\epsilon^5 \right).
\]
This non asymptotic high concentration  expansion of the first moment of the empirical mean on manifolds exhibits an unexpected bias in $1/n$ proportional to the (covariant) gradient of the curvature tensor contracted twice with the covariance matrix. This bias appearing in the small sample size regime was apparently completely unnoticed before.

For the second moment, Theorem \ref{thm:MomentsEmpiricalMean} states 
 that the covariance of the empirical mean is: 
%$\text{Cov}(\bar x_n) = \expect{ \log_{\bar x}(\bar x_n) \otimes \log_{\bar x}(\bar x_n) }$, is:
\[
%\text{Cov}(\bar x_n) = 
  \expect{ \log_{\bar x}(\bar x_n) \otimes \log_{\bar x}(\bar x_n) } =  
  \textstyle \frac{1}{n} \left( \MM_2  -\frac{1}{3}\left(1-\frac{1}{n}\right) 
				\MM_2 \contractcirc ( \contractvarcirc \otimes R(\contractvar, \contractvarcirc)\contractvar + R(\contractvar, \contractvarcirc)\contractvar \otimes \contractvarcirc  )\contract \MM_2   \right)
 + O\left(\epsilon^5\right) .  %% , \textstyle \frac{1}{n^2} \right) .
\]
This non asymptotic high concentration  expansion  of the covariance matrix clearly shows  a modulation of the convergence rate in $1/n$ proportional to the covariance-curvature tensor.  Experiments on constant curvature spaces demonstrate that the expansions are very accurate for variances that are smaller than the curvature, and allows to draw an intuitive understanding of the impact of the curvature on the statistical estimation. 

\paragraph{Link with the BP-CLT}
We compare in Section \ref{sec:BP-CLT} our non-asymptotic expansion to the asymptotic expansion on Riemannian manifolds \cite[Theorem 2.2 p. 1231]{bhattacharya_large_2005}. Rephrased with our notations, the BP-CLT states that, in KKC conditions,  the empirical Fr\'echet mean $\bar x_n$ is a consistent estimator of the population Fr\'echet mean $\bar x$ and the random variable $\sqrt{n} \log_{\bar x} (\bar x_n) \in T_{\bar x} \M$ converges in law to  a normal distribution of mean 0 and covariance $4 \bar H\inv \: \MM_2 \: \bar H\inv$, where the matrix $\bar H$ is the expectation of the Riemannian Hessian of the squared distance $\dist( . , y)^2$ (Theorem \ref{BP-CLT}). In this formula, we see that the expected Hessian is controlling the speed of convergence to the mean.  Thus, the breaking of the BP-CLT in the case of zero eigenvalues of the expected Hessian may also be interpreted as an absence of convergence.

Establishing the high concentration expansion of $\bar H$ shows that both expansions are asymptotically consistent. However, our new non-asymptotic high concentration expansion has additional correction terms for the small data regime while the BP-CLT includes implicit correction terms in the expected Hessian of the squared distance for less concentrated distributions. 
One of the main interests of our new expansion is also to give an intuitive and intelligible interpretation of the expected Hessian of the squared distance  in terms of the curvature and its derivatives.

\paragraph{Modulation of the convergence rate in space forms}
In order to better visualize the influence of the curvature on the empirical mean, we investigate in Section \ref{sec:ModulartionConstantCurvature} the case of isotropic distributions in constant curvature spaces, also called space forms. In this case, both the asymptotic and the high concentration expansions can be reduced to a scalar equation where the rate of convergence of the variance of the empirical mean with respect to the number of samples is modulated by a scalar factor $\alpha = \text{Var}(\bar x_n) {n} / {\sigma^2}$ indicating how much the variance of the empirical mean deviates from the Euclidean case: a modulation factor $\alpha >1$ indicates that the convergence is slower  than in the Euclidean case, while a value $\alpha < 1$ indicates a faster convergence. 

The setup can be further simplified by considering a uniform distribution on a Riemannian hypersphere of radius $\theta$ around the population mean $\bar x$.  For  a large number of samples in a manifold of large dimension, we obtain an  archetypal modulation factor
$\alpha = \frac{\tan^2( \sqrt{ \kappa \theta^2} )}{ \kappa \theta^2}$ on the convergence of the variance of the empirical mean.
We see that the variable controlling the modulation is actually $\kappa \theta^2$, the product of the sectional curvature with the variance. 
For a positive variance-curvature, the modulation of the rate of convergence is larger than one (the convergence is slower) and goes to infinity when $\kappa \theta^2$ approaches $\pi^2 / 4$. This corresponds exactly to the Kendall \& Karcher concentration conditions under which all the results of this paper are restricted. This was expected since a uniform distributions on a Riemannian hypersphere of radius $\pi/2$ on a sphere fails to have a unique mean: the distribution of the empirical Fr\'echet mean converges to a mixture of Diracs rather than concentrating on a point as usual. 
For negative curvature, the modulation factor is below 1, meaning that the convergence is accelerated, and actually goes to zero (an infinite acceleration) for an infinitely negative variance-curvature.
Such a phenomenon has been observed in specific cases for other types of means in negatively curves manifolds but was apparently not recognized so far as a general phenomenon of least-squares in manifolds.  We conjecture that it is related to the phenomenon of stickiness of the Fr\'echet mean in stratified spaces.
 
These theoretical predictions are illustrated by experiments on the 2- and 3-sphere and on the hyperbolic space  of dimension 3 which demonstrate that the formulas that we have obtained are very accurate in their own domains.

\section{Means on Riemannian and affine connection manifolds}

\subsection{Riemannian manifolds}

We consider a differential manifold  $\M$ provided with a smooth scalar product $\scal{.}{.}_{x}$ on each tangent space $T_{x}\M$ at point $x$ of $\M$, called the Riemannian metric. 
In a chart, the metric is fully specified by the dot product of the coordinate vector fields: $g_{ij}(x) = \scal{\partial_i}{\partial_j}$. The Riemannian distance between any two points on $\M$ is the infimum of the length of the curves joining these points. 
Geodesics are defined as the critical points of the energy functional.
% ${\cal E}(\gamma) = \frac{1}{2}\int_0^1 \left\| \partial_{\gamma}\right\|^2\: dt$ 
Geodesics are  parametrized by arc-length in addition to optimizing the length functional. 
%Denoting $[g^{ij}] = [g_{ij}]\inv$ be the inverse of the metric matrix and $\Gamma^i_{jk} = \frac{1}{2} g^{im}\left( \partial_k g_{mj} + \partial_j g_{mk} - \partial_m g_{jk} \right)$ the Christoffel symbols (using Einstein summation convention),  the calculus of variations shows the geodesics satisfy the following second order differential system: \[  \ddot{\gamma}^i + \Gamma^i_{jk} \dot{\gamma}^j \dot{\gamma}^k = 0.\]
% The fundamental theorem of Riemannian geometry states that on any Riemannian manifold there is a unique (torsion-free) connection which is compatible with the metric, called the Levi-Civita (or metric) connection. This connection is determined in a local coordinate system through the  Christoffel symbols:  $\nabla_{\partial_i}\partial_j = \sum_k \Gamma_{ij}^k. \partial_k$. For that choice of connection, shortest path are geodesics ("straight lines"). In the following, we only consider the Levi-Civita connection. 
We assume in this paper that the manifold is geodesically complete, i.e. that the definition domain of all geodesics can be extended to $\R$. This means that the manifold has no boundary nor any singular point that we can reach in a finite time. As an important consequence, the Hopf-Rinow-De~Rham theorem states that there always exists at least one minimizing geodesic between any two points of the manifold (i.e. whose length is the distance between the two points).

From the theory of second order differential equations, 
we know that there exists one and only one geodesic
$\gamma_{(x,v)}(t)$ starting from the point $x$ with the
tangent vector $v \in T_{x}\M$ .
The exponential map at point $x$  maps each tangent vector $v \in
T_{x}\M$ to the point of the manifold that is reached
after a unit time by the geodesic: $ \exp_{x}(v) = \gamma_{(x,v)}(1)$.
The exponential map is locally one-to-one around $0$:  we
denote by $\log_{x}(y)$ its inverse. To shorten formulas, we sometimes use the notation $\lcp{xy}$ instead of $\log_{x}(y)$. 

The cut time $t_{cut}(x,v)$ is the maximal time for which the normal geodesic starting at $x$ with unit tangent vector $v \in T_x\M$ is length minimizing.
By homogeneity, the cut-time can be extended to the tangent bundle $T\M^*$: $t_{cut}(x,v) = t_{cut}(x,v/\|v\|_x) \|v\|_x$, except for null tangent vectors.
The tangent cut-locus is defined as the set of vectors of $T_x\M^* =  T_x\M \setminus \{0\}$ where the distance ceases to be minimizing: $C(x) = \{ t_{cut}(x,v) v | v \in T_x\M^* \}$. The cut-locus is the image of the tangent cut-locus by the exponential map: $\Cut(x) = \exp_x( C(x) )$.
The distance to a point $x$ is $C^2$ except at the cut locus where it is only continuous. 
The tangent cut locus delimits a star shaped domain around $0$ in each tangent space whose interior is called the injectivity domain of the exponential map:
\[
\Inj(x) = \{ tv \: | \: 0 \leq t < t_{cut}(x,v), v \in T_x\M^* \}. 
\]
The injection radius is the infimum of the cut values of the various geodesics emanating from $x$ (i.e. the radius of the largest open geodesic ball included in the injectivity domain, and thus on which $\exp_x$ is a diffeomorphism):
\[
\inj(x) = \min \{ t_{cut}(x,v/\|v\|_x) \: | \: v \in T_x\M^* \}.
\]

There exists multiple notions of convexity in manifolds \cite{berger_panoramic_2003}. We use here the following: an open subset $U \subset \M$ is (strongly) convex if for any points $p,q \in U$, there exists a unique minimal geodesic $\gamma$ joining $p$ and $q$ which belongs entirely to the subset: $\gamma \subset U$. 

%\begin{definition} [Moments of a probability measure] 
%\label{def:moments}
%Let $\mu  \in {Prob}(\M)$ be a probability measure on a Riemannian manifold $\M$. 
%%This probability measure might be uniformly bounded by the Riemannian measure (i.e. such that $\mu(dx)= \mu(x)d\M(x)$) but it may also include mass densities and their derivatives. 
%At the points of $\M$ where there is no mass on the cut locus ($\mu( \Cut(x))=0$), we define the $k$-order moment as the $(k,0)$ tensor:
%\begin{equation}
 %\MM_k\mu(x) = \int_\M \underbrace{\log_x(y) \otimes \log_x(y) \ldots  \otimes  \log_x(y)}_{\text{$k$ times}} \: \mu(dy)
%= \int_\M \log_x(y)^{\otimes k} \: \mu(dy)
%\end{equation}
%\end{definition}
%In general, the tensor field is non-smooth at points where the cut locus has a non-zero mass.
%If the density $\mu$ is uniformly bounded by the Riemannian measure (i.e. if $\mu(dx)= \mu(x)d\M(x)$), then this integral defines a smooth $(k,0)$ tensor field over $\M$ since the cut locus of the point $x$ has null measure. The 0-th order moment $\MM_0\mu = \int_{\M} d\mu(x) =1$ is unit by definition. The first order moment $\MM_1\mu(x) =  \int_{\M}\lcp{xx_i}\: d\mu(x)$ is a vector field on the manifold $\M$.

\subsection{Convex affine manifolds}

An affine manifold is a differential manifold $\M$ endowed with an affine connection $\nabla$.   Important classes of affine manifolds are: Riemannian manifolds with their Levi-Civita connection, Lie groups with their canonical symmetric space structure, 
%thus endowed with the symmetric Cartan-Schouten bi-invariant connection, 
and more generally affine symmetric spaces. The connection allows us to define geodesics as auto-parallel curves, or zero acceleration curves. These geodesics have an affine parametrization that  measures relative distances along each geodesic, but there is no reference length to compare relative distances along different geodesics.  Convexity may be defined as in Riemannian manifolds:  an open set $\U \subset \M$ is convex if for every pair of points $x,y \in \U$ there exists a unique $\nabla$-geodesic joining $x$ and $y$ lying entirely in $\U$ and that depends smoothly on its endpoints. This allows us to define the logarithm in a unique way within this subset. 
In that case, $(\U, \nabla)$ is called a convex (sub)-manifold. Whitehead Theorem  tells us that there exists a convex neighborhood at each point of an affine manifold.

\subsection{Curvature of a connection}

In order to derive  Taylor expansions in manifolds, we  need to specify the notations for torsion and curvature operators and their coordinates in charts. The torsion tensor $T(X,Y) = \nabla_X Y -\nabla_Y X - [X,Y] = -T(Y,X)$  measures how the skew-symmetric part of the connection differ from the Lie derivative ${\cal L}_X Y = [X,Y]$. The connection is torsion free if the torsion tensor vanishes identically. 
Two connections have the same geodesics if they have the same symmetric part $(\nabla_X Y + \nabla_Y X)/2$. i.e. if they only differ by torsion. Because means and barycenters in manifolds only involve geodesics and not parallel transport, we can restrict our attention to torsion free (also called symmetric) connections. This notably simplifies the expression of the Taylor expansions. In the sequel, all connections are assumed to be symmetric.

\paragraph{The Riemann curvature tensor in affine connection manifolds}
The curvature of an affine manifold is described by the $(1,3)$ curvature tensor $R:T\M\times T\M \times T\M\to TM$. It is defined from the covariant derivative by evaluation on vector fields $u,v,w$:
\begin{equation}
  R(u,v)w  = \nabla_u\nabla_v w - \nabla_v \nabla_u w -\nabla_{[u,v]} w.
  \label{eq:CurvatureTensor}
\end{equation}
The curvature can be interpreted as the difference between parallel transporting the vector $w$ along an infinitesimal parallelogram with sides given by  $u$ and $v$. 
Two different sign conventions exist for the curvature tensor: the above definition  is the one used in a number of reference books in physics and mathematics \cite{misner_gravitation_1973,lee_riemannian_1997,postnikov_geometry_2010}. %,jost_riemannian_2011}. %klingenberg82
% sternberg
Other authors use a minus sign to simplify some of the tensor notations  \cite{spivak_differential_1979,oneill_semi-riemannian_1983,carmo_riemannian_1992,berger_panoramic_2003}. %Gallot:93
There exists moreover different conventions for the order of the tensors subscripts. We use here the convention $[R(u,v)w ]^a = R^a_{bcd} u^c v^d w^b$ so that the tensor can be written using the Christoffel symbols $\Gamma^a_{bc}$: 
%and/or a minus sign in the sectional curvature defined below (see e.g. the discussion in \cite[p. 399]{gallier_notes_2016}). 
\begin{equation}
R^a_{bcd} = dx^a ( R(\partial_c, \partial_d) \partial_b) 
= \partial_c \Gamma^a_{db} -\partial_d \Gamma^a_{cb} + \Gamma^a_{ce}\Gamma^e_{db} - \Gamma^a_{de}\Gamma^e_{cb}, 
\label{eq:CurvatureCoord}
\end{equation}
where $\partial_i = \partial/\partial x_i$ are coordinate vector fields. 
 
The curvature tensor has several symmetries that we will use to simplify expressions: it is skew-symmetric in the first two variables, and we we have the First and second Bianchi identities. 
\begin{eqnarray*}
 R(u,v)=-R(v,u) & \text{or} & R^a_{bcd} = - R^a_{bdc} ; \\
 R(u,v)w + R(v,w)u + R(w,u)v=0 & \text{or} & R^a_{[bcd]} = R^a_{bcd} + R^a_{cdb} + R^a_{dbc} =  0 ; \\
 (\nabla_u R)(v,w) + (\nabla_v R)(w,u) + (\nabla_w R)(u,v)=0 & \text{or} & 
R^a_{b[cd;e]} = \nabla_e R^a_{bcd} + \nabla_c R^a_{bde} + \nabla_d R^a_{bec} =  0.
%R^a_{b[cd;e]} = R^a_{bcd;e} + R^a_{bde;c} + R^a_{bec;d} =  0.
\end{eqnarray*}

\paragraph{Additional symmetries in a Riemannian manifold}
In a Riemannian manifold, we can lower the first coordinate with the metric to obtain 
the $(0,4)$ version of the Riemannian curvature tensor:
\begin{equation}
R(u,v,w,z) = \scal{R(w,z)v}{u} 
\quad\text{or in coordinates} \quad
R_{abcd} = g_{ae} R^e_{bcd}.
\label{eq:CovariantCurvatureTensor}
\end{equation}
This tensor inherits the above symmetries. Additionally, it is symmetric in the first and last two variables, which implies being skew-symmetric in the last two variables:
\begin{eqnarray*}
 \scal{R(u,v)w}{z} = -\scal{R(u,v)z}{w} & \text{or}  & R_{abcd} = R_{cdab} ; \\
\scal{R(u,v)w}{z} = \scal{R(w,z)u}{w} & \text{or} & R_{abcd} = -R_{abdc} .
\end{eqnarray*}
This means in particular that its contraction with any symmetric tensor in the first two or last two indices is zero.

Note that in an affine manifold with an auxiliary metric $g$ we can also lower the first index of the curvature tensor, but this metric should be kept in all computations as it is not the identity in a normal coordinate system and its covariant derivative is not zero as with the Levi-Civita connection.

\paragraph{Sectional curvature}

In Riemannian manifolds, the sectional curvature $\kappa(u,v)(x)$ measures the Gaussian curvature (the product of the principal curvatures) in the 2-planes of $T_x\M$ generated by the vectors $u$ and $v$.  It can be expressed from the curvature tensor by: 
\begin{equation}
    \kappa(u,v)(x) = \frac{\scal{R(u,v)v}{u}_x}{\|u\|_x^2\|v\|_x^2 - \scal{u}{v}_x^2}
    \, .
  \label{eq:sec_curvature}
\end{equation}
The upper bound of the sectional curvature plays a very important role for convexity in Riemannian manifolds, as we will see below.

\subsection{Moments of a probability measure}

\begin{definition} 
\label{def:moments}
Let $\mu  \in {Prob}(\M)$ be a probability measure on a  manifold $\M$. 
%This probability measure might be uniformly bounded by the Riemannian measure (i.e. such that $\mu(dx)= \mu(x)d\M(x)$) but it may also include mass densities and their derivatives. 
%At the points of $\M$ where there is no mass on the cut locus ($\mu( \Cut(x))=0$) for a Riemannian manifold, and within ??? for a convex affine connction manifold, we define 
The $k$-order moment of $\mu$ is the $(k,0)$ tensor:
\begin{equation}
 \MM_k\mu(x) = \int_\M \underbrace{\log_x(y) \otimes \log_x(y) \ldots  \otimes  \log_x(y)}_{\text{$k$ times}} \: \mu(dy)
= \int_\M \log_x(y)^{\otimes k} \: \mu(dy)
\end{equation}
\end{definition}
For a general distribution in a Riemannian manifold, the tensor field is not defined (and non-smooth) at points where the cut locus has a non-zero mass.
If the density $\mu$ is uniformly bounded by the Riemannian measure (i.e. if $\mu(dx)= \mu(x)d\M(x)$), then this integral defines a smooth $(k,0)$ tensor field over $\M$ since the cut locus of the point $x$ has null measure. In an affine manifold,  the support of $\mu$ needs to be limited to an convex neighborhood so that the logarithm is well defined and all moments are smooth within this neighborhood.
The 0-th order moment $\MM_0\mu = \int_{\M} d\mu(x) =1$ is unit by definition. The first order moment $\MM_1\mu(x) =  \int_{\M}\lcp{xx_i}\: d\mu(x)$ is a vector field on the manifold $\M$.

In the following, we only consider distributions that have support in a convex affine manifold $(\M, \nabla)$ that is fixed once for all. Such a convex manifold is diffeomorphic to an open set of $\R^d$. 
 Cartan-Hadamard manifolds (complete, simply connected manifolds with sectional curvature less than or equal to 0) are classical convex complete Riemannian manifolds. Examples of incomplete convex manifolds are given by small geodesic balls in Riemannian manifolds or small balls centered at the origin in an exponential chart in an affine manifold \cite{arnaudon_barycenters_2005}.

\subsection{Means on Riemannian manifolds}

The classical mean value of random numeric values with distribution $\mu$ is defined through an integral $\bar x = \int x \: \mu(dx)$. Because $\mu$ is normalized, this definition can be rewritten as an implicit barycentric equation: $ \int (x -\bar x) \: \mu(dx)$. With this formulation, it is obvious that this notion is affine and not metric. In the context of probability and statistics, Maurice Fr\'echet was the first to unify several type of typical statistical values, in particular the mean and the median, and to  generalize them to abstract spaces. In a preparatory work, he first investigated  \cite{frechet_valeurs_1943} different ways to compute mean values of random triangles, independently of their position and orientation in space, with experiments to confront the theoretical results to real world data. In this respect, he may be considered as pioneering the statistical study of shapes. In a second study, motivated by the study of random curves, he first introduced a mean value and a law of large numbers defined by a generalization of the integral to normed vector (Wiener or Banach) space.  
Finally, Fr\'echet considered in \cite[p.233]{frechet_les_1948} a family of central values (including the mean and the median) and its generalization to random elements in abstract metric spaces: 
\begin{definition}[Fr\'echet mean in a metric space \cite{frechet_les_1948}] %$ $ \\
The $p$-mean (typical position of order $p$ according to Fr\'echet) of a distribution  (a random element) $\mu$ in an abstract metric space $\M$ is set of minima of the $p$-standard deviation : 
\begin{equation} \textstyle
%% \expect{x_i} = 
\text{p-Mean} (\mu) =  \left\{ \underset{y \in \M}{\arg \min} \; \int_\M \dist(x, y)^p \: \mu(dx) \right\}. 
\end{equation}
The case $p=2$ corresponds in vector spaces to the arithmetic mean, the case $p=1$  to the median (``equiprobable value'' in Fr\'echet's words).
\end{definition}
The first key contribution of Fr\'echet was to consider many different types of typical elements, including of course the mean but also the median. Fr\'echet considered mainly the case $p\geq 1$, but he observed that many of the properties could be also generalized to $0<p<1$. His second revolutionary idea was to considered  a set of mean elements rather than one unique mean. This idea was later  developed by Ziezold \cite{ziezold_expected_1977} with a strong law of large numbers for sets of random elements in separable finite quasi-metric spaces. These two innovations justify the name of Fr\'echet mean that is used in geometric statistics.

For smooth differential geometric spaces like Riemannian manifolds, and restricting to the classical 2-mean, Berger \cite[p.235]{berger_panoramic_2003} reported that ``the existence of a unique center of mass in the large for manifolds with non-positive curvature was proven and used by \'Elie Cartan back in the 1920's.''
In order to find the fixed point of a group of isometries, Cartan indeed showed in \cite{cartan_cons_1928}  that the sum of the square distances from one point to a finite number of points has a unique minimum in simply connected Riemannian manifolds with non-positive curvature\footnote{Note III on normal spaces with negative or null Riemannian curvature, p. 267.} (now called a Hadamard space). This result was extended in \cite{cartan_groupes_1929} to closed subgroups of isometries. It is obvious in this text that Cartan is only using the uniqueness of the minimum of the sum of square distances as a tool in the specific case of negative curvature Riemannian manifolds and not as a general definition of the mean on manifolds as usually thought in probability or statistics.

For similar purposes, Grove and Karcher extended this idea in 1973 to positive curvature manifolds  for distributions with sufficiently small support, typically within convex balls so that the mean exists and is unique \cite{grove_how_1973}. The notion was coined {\em the Riemannian center of mass}. In this publication and in successive ones, Karcher and colleagues determined Jacobi field estimates to find optimal conditions for the convexity of the ball that support this definition. The Riemannian barycenter is commonly refered to as describes in  \cite{karcher_riemannian_1977} but the most complete description of the related properties is certainly found in \cite{buser_gromovs_1981}, where a notion of barycenter in affine connection manifolds is also worked out. A good historical note on the history of the Riemannian barycenter is given in \cite{afsari_riemannian_2011} and by Karcher himself in \cite{karcher_riemannian_2014}. 

In all these works, the Riemannian center of mass is unique by definition. Considering a set-valued barycenter on an affine connection manifold was the contribution of Emery and Mokobodzki:
\begin{definition}[Exponential barycenters in an affine connection manifold \cite{emery_sur_1991}]
We consider a probability measure $\mu$ with support in a convex local neighborhood $\U$ of an affine connection manifold $(\M, \nabla)$.  Exponential barycenters of the probability measure $\mu$ are the points 
\[
\text{Bary}_{\U}(\mu) = \left\{ x \in \U \: | \: \MM_1(x) = \int_{\M} \log_x(y) \: \mu(dy) =0 \right\}.
\] 
\end{definition}
The non-uniqueness of the expectation of a random variable considerably extended the usability of this definition, in particular in positive curvature manifolds. We notice that the notion remains purely affine, provided that the distribution has a support on a convex neighborhood in which the logarithm can be defined uniquely. In the Riemannian case, \cite{emery_sur_1991} showed that exponential barycenters are the critical points of the variance. Thus exponential barycenters contain in particular the minimum of the variance $\sigma^2(x) = \int_\M \dist^2(x,y)\: \mu(dy)$ (sometimes called the Fr\'echet functional), except if the variance is not differentiable at this point. 

The optimal conditions for the existence and uniqueness of the minimum of the variance have been studied in \cite{karcher_riemannian_1977,buser_gromovs_1981,kendall_probability_1990,le_locating_2001,le_estimation_2004}. The result has been extended to Fr\'echet $p$-means defined as the minima of the $p$-variance in \cite{afsari_riemannian_2011,yang_riemannian_2010}. The Hessian of the squared distance to a point plays a key role in these conditions. Karcher looked for a ball where this Hessian is positive definite everywhere, while Kendall relaxed this condition to the locus of the minimum. We will see in Section \ref{sec:ArchetypalModulation} that the inverse of the expected Hessian also controls the speed of convergence to the mean, so that the non-uniqueness may also be interpreted as an absence of convergence.  
\begin{theorem}[Karcher \& Kendall Concentration (KKC) Conditions] 
%$ $\\
\label{KKCC}
Let $\M$ be a geodesically complete Riemannian manifold with injection radius $\inj(x)$
Let $\mu$ is a probability distribution on $\M$ whose support is contained in a closed regular geodesic ball $\bar B(x,r)$ of radius $r < \frac{1}{2} \inj(x)$. We assume moreover that the upper bound $\kappa = \sup_{y\in B(x,r), u \not = v \in T_y\M} \kappa(u,v)(y) $ of the sectional curvatures in the ball satisfies $\kappa < \pi^2 /(2r)^2$. This second condition is always verified on spaces of negative curvature and it specifies a maximal radius $r^* =  \frac{\pi}{2\sqrt{\kappa}}$ when there is positive sectional curvature.  
These concentration assumptions ensure that the variance has a unique global minimum that belongs to the ball $\bar B(x,r)$.
\end{theorem}
 
In order to differentiate the different notions of means in Riemannian manifolds, it is usual in geometric statistics to name Fr\'echet mean the set of global minima of the variance, Karcher mean the set of local minima, and exponential barycenters the set of critical points satisfying the implicit equation $\MM_1(x) =0$. 
It is clear that all these definition boil down to the same unique point within the ball $B(x,r)$ in KKC conditions for the classical 2-mean, although some local minima and critical points may exist outside this ball.
 Throughout this paper, we assume that the support of all the distributions that we consider on Riemannian manifolds are included in a regular geodesic ball $B(x_0,r)$ of diameter $2 r<\varepsilon$ that satisfies the KKC conditions. As a consequence, one can conclude that there is a unique exponential barycenter $\bar x$ included in the ball, necessarily the Fr\'echet mean \cite{karcher_riemannian_1977,kendall_probability_1990}. The maximal diameter $\varepsilon$ of the support of the distribution is used in this paper as the scale variable to control Taylor expansions.

\subsection{Means in convex affine connection manifolds}

To define the mean of a probability distribution $\mu$ in a convex manifold $(\M, \nabla)$,  we cannot rely on the Fr\'echet mean since there is no distance. However,  the notion of exponential barycenter still makes sense because the affine logarithm is well defined \cite{emery_sur_1991}: they are  the zeros of the first moment field ${\mathfrak M}_1 (\mu) = \int_{\M} \log_{x}(y) \: \mu(dy)$ for $x\in M$.  
%We notice that there is an implicit dependency of the convex set $\U$ considered in that definition. 
This definition, studied in \cite{buser_gromovs_1981,emery_sur_1991}, is close to  the Riemannian center of mass but uses the logarithm of the affine connection in a convex domain instead of the Riemannian logarithm in the injectivity domain.

Distributions with compact support in a convex affine manifold have at least one exponential barycenter. Moreover, exponential barycenter are stable by affine diffeomorphisms (connection preserving maps, which thus preserve the geodesics and the normal convex neighborhoods) \cite{buser_gromovs_1981}. However, the classical notion of convexity is not sufficient to ensure the uniqueness of the mean. For that purpose,  stronger convexity conditions have been proposed.

\begin{definition}[Arnaudon \& Li convexity (ALC) conditions \cite{arnaudon_barycenters_2005}]
Let $(\M, \nabla)$ be an affine manifold.
A separating function on $\M$ is a convex function $\rho: \M \times \M \rightarrow \R^+$ vanishing exactly on the diagonal of the product manifold (considered as an affine manifold with the direct product connection). Here, convex means that the restriction of $\rho(\gamma(t))$ to any geodesic $\gamma(t)$ of $\M \times \M$ is a convex function from $\R$ to $\R^+$. 
A manifold which carries a smooth separating function $\rho$ such that
\[ c \: \dist(x,y)^p \leq \rho(x,y) \leq C \: \dist(x,y)^p, \]
for some constants $0 < c < C$, some positive integer $p \geq 2$ and some auxiliary Riemannian distance function $\dist$ is called a manifold with $p$-convex geometry.
On a manifold with $p$-convex geometry, every compactly supported probability measure has a unique exponential barycenter. 
\end{definition}
%A separating function generalizes the separating property of a distance ($\rho(x,y)=0 \Leftrightarrow x=y$) and the ever increasing value when $y$ gets away from $x$ along geodesics but does not respect the triangular inequality. 
Whitehead theorem tells us in essence that any point in an affine connection manifold has a convex neighborhood with 2-convex geometry. In Riemannian manifolds, geodesic balls in KKC conditions have 2-convex geometry. These conditions will be sufficient for our goal of studying the behavior of the empirical means with high concentration. We should note that the notion is too strong in general since there are examples of manifolds where one can prove uniqueness although they do not have $p$-convex geometry for any $p$. This has motivated the definition of CSLCG (convex, with semi-local convex geometry) manifolds in \cite{arnaudon_barycenters_2005} but we will not need it in this paper.

Throughout this paper, we only consider distributions that satisfy the ALC condition, i.e. compact distributions in a convex affine manifold $(\M, \nabla)$ with $p$-convex geometry. With this extra structure, we have an auxiliary Riemannian metric on $\M$ whose norm in $T_x\M$ is denoted $\|.\|_x$  and whose Riemannian distance is denoted $\dist$. Thanks to this distance, we can compute the diameter  $\varepsilon$ of the support of the distribution, which is then used as in Riemannian manifolds as the typical scale variable to control Taylor expansions.

\section{Taylor expansions in Riemannian and affine manifolds \label{sec:Taylor}}

% \subsection{Taylor expansions in manifolds}
In order to analyze the variance of a random element in a Riemannian manifold,  we may compute the Taylor expansion of the squared Riemannian distance functions in a local coordinate system. This requires to have an expression for the geodesics. 
 In a normal coordinate system centered at $x$, a geodesic starting at $x$ with tangent vector $v$ is a straight line: $\log_x( \exp_x(tv)) = tv$. However, the geodesic starting at $x_v = \exp_x(v)$ and ending at $x_w = \exp_x(w)$ deviates from the straight line from $v$ to $w$ in our chart because of curvature. Geodesic may be approximated by a polynomial expansion in the variables $v$ and $w$  that solves the geodesic equation. 

The expansion of the Riemannian logarithm can be based on the Taylor expansion of the Riemannian metric in a normal coordinate system. Such an extension is well known for the first orders, but going to orders higher than 4 is computationally much more involved and requires computational algebraic methods such as the ones developed by Brewin \cite{brewin_riemann_2009}. Part of the results presented in this paper were originally developed in 2015 using this method. However, the intelligibility of the formulas involving lots of terms and indices was difficult.

In this paper, we base our Taylor expansions on a more geometric formulation based on a coordinate free expansion of the composition of two exponential maps developed by Gavrilov \cite{gavrilov_algebraic_2006,gavrilov_double_2007}. This formula only involves the curvature and torsion tensors and their covariant derivatives at the development point. Moreover, the proposed extension relies only on the connection and not on the Riemannian metric, which makes it immediately suitable for the general affine connection case.

\subsection{Gavrilov's expansions of the double exponential in affine manifolds}

\begin{figure}[tb!]
	\includegraphics[width=7cm]{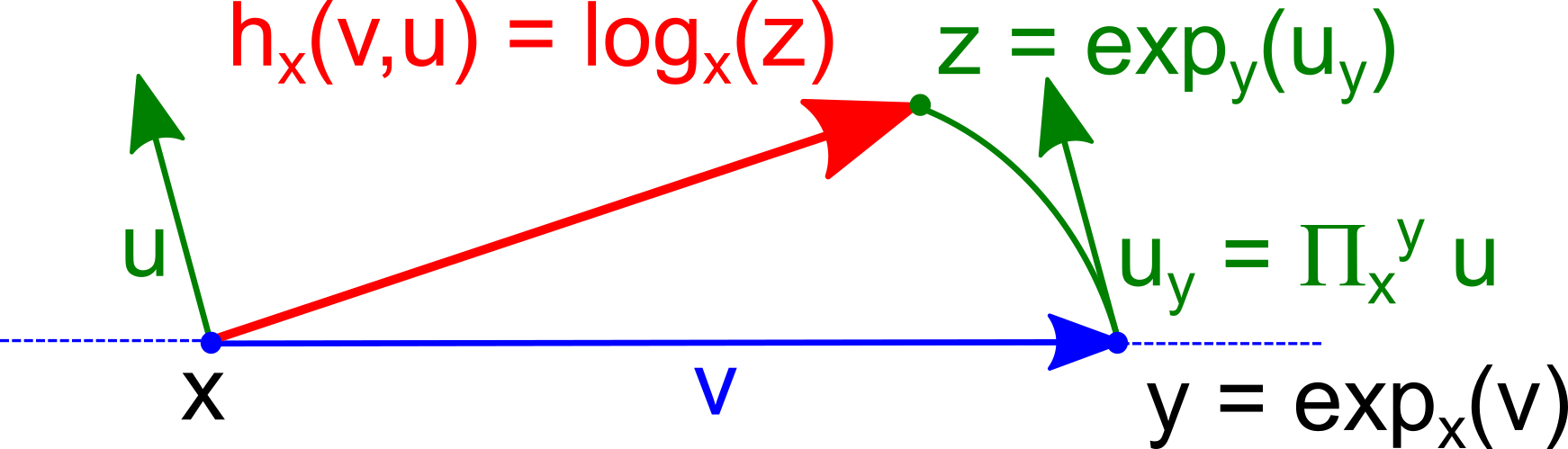}\hfill
	\includegraphics[width=7cm]{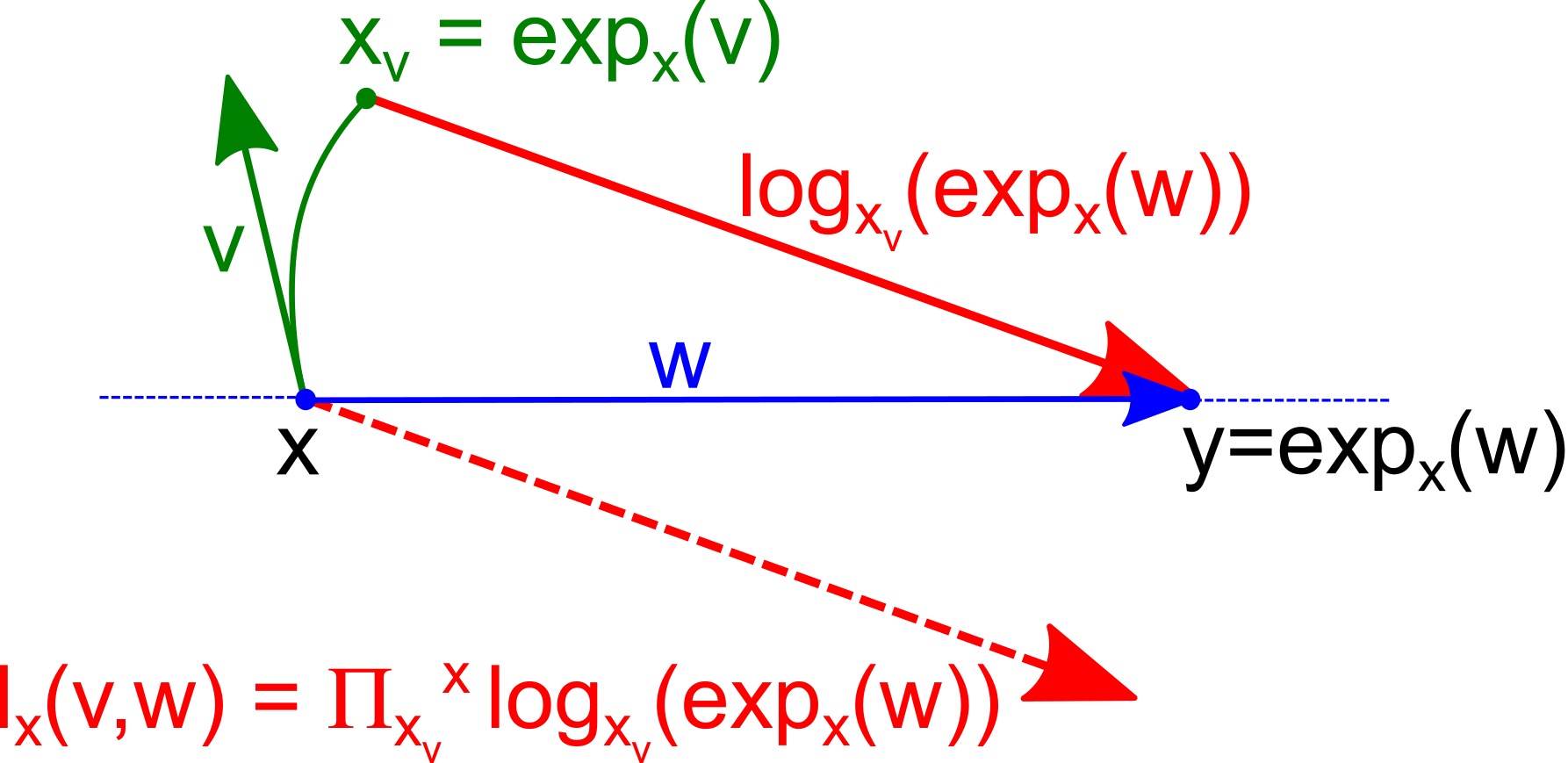}
\caption{{\bf Left:} the log of the composition of two exponentials (BCH-type formula) in a normal coordinate system at x.
{\bf Right:} The neighboring log of a fixed point $\exp_x(w)$ when the foot-point $x$ is moved along the geodesic $x_v= \exp_x(v)$, in a normal coordinate system at x.
}\label{Fig:BCH}
\end{figure}

In the spirit of the Baker-Campbell-Hausdorff (BCH) formula for Lie groups, Gavrilov \cite{gavrilov_algebraic_2006,gavrilov_double_2007} developed a coordinate free expansion of the composition of two exponential maps. 
The double exponential $\exp_x(v,u) = \exp_{\exp_x(v)}( \Pi_x^{\exp_x(v)} u)$ corresponds to a first geodesic shooting from the point $x$ along the vector $v$, followed by a second geodesic shooting from $y= \exp_x(v)$ along the parallel transport $\Pi_x^y u$ of the vector $u$ along the first geodesic (Fig. \ref{Fig:BCH}, left). This expansion holds in general affine connection manifolds, and has a surprisingly simple coordinate free formulation. We formulate it here in the torsion free case.

%\begin{figure}[tb!]
%\centering
	%\includegraphics[width=8cm]{Figures/BCH.png}
%\caption{The log of the composition of two exponentials (BCH-type formula) in a normal coordinate system at x.}\label{Fig:BCH}
%\end{figure}

\begin{theorem}[The double exponential expansion \cite{gavrilov_algebraic_2006,gavrilov_double_2007}] \label{thm:DoubleExpExpansion}
In a torsion-free affine connection manifold, the log of the double exponential $h_x(v,u) = \log_x(\exp_x(v,u))$ has the following series expansion in $v$ and $w$ at order 5:
\begin{equation}\label{Eq:DoubleExp}
\begin{split}
h_x(v, u) =\,  &v + u + \frac{1}{6}R(u,v)v + \frac{1}{3}R(u,v)u 
+ \frac{1}{24}(\nabla_v R)(u,v) (2v +5u)  \\ &+ \frac{1}{24}(\nabla_u R)(u,v) (v + 2u)
+O(5),
\end{split}
\end{equation}
where  $O(5)$ represents polynomial terms of order 5 or more in u and v.  
\end{theorem}
In the right-hand side, the tensor values are taken at $x$. By extending the vectors $u$ and $v$ of $T_x\M$ to vector fields in the neighborhood of $x$ using parallel transport, we see that the series is valid point-wise at any point in that neighborhood.  

The fundamental idea of this expansion is the use of parallel transport to identify tangent spaces rather than the differential of the exponential as is implicitly done in Taylor expansions in normal coordinate systems \cite{brewin_riemann_2009}. This trick drastically simplifies the expressions with coordinate free expressions. We notice that $h_x$ is a mapping of vector spaces from $T_x\M \times T_x\M$ to $T_x\M$.

\subsection{Moving the foot-point: the neighboring log expansion}
The double exponential formula may be used to evaluate how the log of a fixed point $y=\exp_x(w)$ is modified when $x$ is moved along the geodesic $x_v= \exp_x(v)$: we define the neighboring log as $l_x(v, w) = \Pi_{x_v}^x \log_{x_v}(\exp_x(w))$ (Fig. \ref{Fig:BCH}, right). %\ref{Fig:NLog}). 
%\begin{figure}[htb!]
%\centering
	%\includegraphics[width=8cm]{Figures/Nlog.png}
%\caption{The neighboring log of a fixed point $\exp_x(w)$ when the foot-point $x$ is moved along the geodesic $x_v= \exp_x(v)$, in a normal coordinate system at x.}\label{Fig:NLog}
%\end{figure}
This amounts to say that the vector $u=l_x(v, w)$ is solution of 
%$\Pi_x^{x_v} u = \log_{x_v}(\exp_x(w))$ or after taking the exponential: 
$\exp_{x_v}(\Pi_x^{x_v} u) = \exp_x(w)$. We recognize the double exponential 
$\exp_x(v,u)$ on the left-hand side.
%$\exp_x(v,u) = \exp_{\exp_x(v)}( \Pi_x^{\exp_x(v)} u)$. 
Thus, determining $u=l_x(v, w)$ amounts to solve  $h_x(v,u) = \log_x(\exp_x(v,u)) = w$. 

% solving h_x(v,w)=u
The first order solution of $h_x(v,u)=w$ is obviously $u=w-v + O(2)$. To find the second order, we assume that $u= w-v + u_2 + O(3)$ where $u_2$ is an unknown polynomial expression of order 2 in $v$ and $w$. Plugging this value in Eq.\ref{Eq:DoubleExp} gives $u_2=O(3)$ since the curvature terms are of order 3. Assuming a third order term $u_3$, we now get: 
\[
h_x(v, w-v+u_3 +O(4)) %= w+u_3 + \frac{1}{6}R(w ,v)v + \frac{1}{3}R(w,v)(w-v) + O(4)
= w+u_3 - \frac{1}{6}R(w ,v)v + \frac{1}{3}R(w,v)w + O(4).
\]
Thus we find $u_3 = \frac{1}{6}R(w ,v)(v-2 w)$. Finally, assuming a forth-order term $u_4$ in $u= w-v + \frac{1}{6}R(w ,v)(v-2 w) + u_4 +O(5)$ gives:
\[
\begin{split}
h_x(v, u)  
% = w + \frac{1}{6}R(w ,v)(v-2 w) + u_4 
% + \frac{1}{6}R(w,v)v + \frac{1}{3}R(w,v) (w-v)
% + \frac{1}{24}(\nabla_v R)(w,v) (5w-3v) + \frac{1}{24}(\nabla_{(w-v)} R)(w,v) (2w-v)
% +O(5)
= w  + u_4 
+ \frac{1}{24}(\nabla_v R)(w,v) (5w-3v) + \frac{1}{24}(\nabla_{(w-v)} R)(w,v) (2w-v)
+O(5),
% &= w  + u_4 + \frac{1}{24}(\nabla_v R)(w,v) (3w-2v) + \frac{1}{24}(\nabla_{w} R)(w,v) (2w-v) +O(5)
\end{split}
\]
so that $h_x(v, u)=w$ implies  $u_4 = \frac{1}{24}(\nabla_v R)(w,v) (2v-3w) + \frac{1}{24}(\nabla_{w} R)(w,v) (v-2w)$. 

\begin{theorem}[The neighboring log expansion] \label{thm:NLogExpansion}
In a torsion-free affine connection manifold, the neighboring log  $l_x(v, w) = \Pi_{x_v}^x \log_{x_v}(\exp_x(w))$ has the following series expansion in $v$ and $w$ at order 5:
\begin{equation}
\begin{split}
l_x(v, w) = &  w-v + \frac{1}{6}R(w ,v)(v-2 w)  \\
&+ \frac{1}{24}(\nabla_v R)(w,v) (2v-3w) 
 + \frac{1}{24}(\nabla_{w} R)(w,v) (v-2w) +O(5) .
\label{eq:NLogExpansion}
\end{split}
\end{equation} 
\end{theorem}

\subsection{Expansion of the Riemannian distance}
With the neighboring log expansion, we can now come back to the expansion of the squared Riemannian distance between two points $x_v=\exp_x(v)$ and $x_w= \exp_x(w)$ that are close to $x$:
\[
\begin{split}
\dist(x_v, x_w)^2 & = \| \log_{x_v}(x_w) \|^2_{x_v} = \| \Pi_{x_v}^x \log_{x_v}(x_w) \|^2_{x} = \| l_x(v, w) \|_x^2 \\
& = \|w-v\|^2_x +\frac{1}{3} \scal{R(w ,v)(v-2 w)}{w-v}_x \\
& \quad + \frac{1}{12}\scal{(\nabla_v R)(w,v) (2v-3w)}{w-v}_x \\
& \quad + \frac{1}{12}\scal{(\nabla_{w} R)(w,v) (v-2w)}{w-v}_x 
+O(6)
.
\end{split}
\]
We now exploit the skew symmetry  
$\scal{R(u,v)w}{z} = -\scal{R(u,v)z}{w}$ to simplify the terms of order 4:
%and the first Bianchi identity $R(u,v)w + R(v,w)u + R(w,u)v)=0$:
$
\scal{R(w ,v)(v-2 w)}{w-v} = 
\scal{R(w ,v)v}{w} 
% + \scal{R(w ,v)(v)}{v} -2 \scal{R(w ,v)(w)}{w}
+2 \scal{R(w ,v)w}{v} 
= \scal{R(w ,v)w}{v} .
$
The derived skew-symmetry for the covariant derivative $\scal{(\nabla_v R)(w,v) u}{z} =- \scal{(\nabla_v R)(w,v) z}{u}$ allows us to simplify the terms of order 5:
\[\textstyle
\scal{(\nabla_v R)(w,v) (2v-3w)}{w-v}_x = 
 % \scal{\nabla_v R(w,v) (2v-3w)}{w}_x - \scal{\nabla_v R(w,v) (2v-3w)}{v}_x =
2\scal{\nabla_v R(w,v) v}{w}_x +3 \scal{\nabla_v R(w,v) w}{v}_x = \scal{\nabla_v R(w,v) w}{v}_x .
\]
\[\textstyle
\scal{(\nabla_{w} R)(w,v) (v-2w)}{w-v}_x  = 
% \scal{\nabla_{w} R(w,v) (v-2w)}{w}_x - \scal{\nabla_{w} R(w,v) (v-2w)}{v}_x =
\scal{\nabla_{w} R(w,v) v}{w}_x +2 \scal{\nabla_{w} R(w,v) w}{v}_x 
= \scal{\nabla_{w} R(w,v) w}{v}_x .
\]
Finally, we obtain the Taylor expansion of the squared geodesic distance in a Riemannian manifold between two points $x_v = \exp_x(v)$ and $x_w = \exp_x(w)$ that are close to $x$:
\begin{equation}
\label{eq:TaylorSquareDist}
\begin{split}
\dist^2(x_v, x_w) = \, & \| w - v \|^2_x  + \frac{1}{3} \left< R(w,v)w \, , \, v \right>_x
+\frac{1}{12}  \left< \nabla_{(v+w)} R(w,v)w\, , \, v\right>_x  
+O(6)  .
\end{split}
\end{equation}

In coordinates, using the orthonormal property of a normal coordinate system, this reads: $\dist^2(x_v, x_w) =  (w - v)^a   (w - v)_a  + \frac{1}{3} R_{abcd} v^a w^b w^c v^d 
+\frac{1}{12} \nabla_e R_{abcd} v^a w^b w^c v^d (w+v)^e  +O(6)$.

\section{Expansion of the first moment of the mean \label{sec:FrechetMeanExpansion}}

Let $\mu$ be a probability distribution with support in a convex neighborhood $\U$ of diameter less than $\varepsilon$ respecting the KKC conditions in a Riemannian manifold or the ALC conditions in a convex affine manifold. In such conditions, all the definitions of the mean boils down to the unique zero of the first moment: $\MM_1\mu(x) =  \int_{\U}\log_x(y)\: \mu(dy)=0$. The goal is this section is to find a polynomial approximation of the solution of this equation in the neighborhood of a point $x \in \U$.

For that purpose, we parametrize the points of $\U$ by $x_v = \exp_x(v)$ with $v \in V_x = \log_x(\U) \subset T_x\M$. Since the tangent mean $\MM_1\mu(x_v)$ is a vector field on $\U$, its value belongs to a different tangent space $T_{x_v}\M$ for each value of $v \in V_x$.  Thus, establishing a Taylor expansion requires to express them all in a single chart. A natural choice is a normal coordinate system at $x$, but the resulting expansion will still depend on that specific coordinate system.  In order to have a unique image space, we may also parallel transport each tangent vector from $x_v$ to $x$ to obtain the following mapping of vector space.

\begin{definition}[Recentered tangent mean map] 
\label{def:RecenterdMeanMap}
Let $\mu$ be a probability distribution with support in a convex open subset $\U$ of an affine manifold $\M$. Then, for any vector $v \in V_x = \log_x(\U) \subset T_x\M$ for which $\MM_1\mu(\exp_x(v))$ is finite, we define the recentered tangent mean map as:
\begin{equation}
\NN_x^\mu(v) = \Pi_{\exp_x(v)}^x \MM_1\mu(\exp_x(v)) \in  T_x\M.
\label{eq:N}
\end{equation} 
This is a mapping of vector spaces from $V_x \subset T_x\M$ to $T_x\M$ whose zeros parametrize the exponential barycenters of $\mu$. 
\end{definition}
The mapping can be extended to the whole injection domain $\Inj(x)$ in the Riemannian case. 
The zeros of the recentered tangent map parametrize the exponential barycenters of $\mu$ because parallel transport is an isometry of tangent spaces in Riemannian manifolds: $\NN_x^{\mu}(v)=0 $ if and only if $\MM_1\mu(\exp_x(v))=0$. This property also holds for affine manifolds since parallel transport is an isomorphism of tangent spaces in that case. In KKC or ALC conditions, we know that the zero is unique.

Since we only deal with one unique probability distribution $\mu$ in the sequel, we simplify the notations below by dropping the dependency on $\mu$ in $\MM_k$ and $\NN_x$.

\subsection{Taylor expansion of the recentered mean map}
Using the neighboring log formula (Eq. \ref{eq:NLogExpansion}), we can compute the series expansion of the recentered mean map. By linearity of the integral, we have:
\[
\begin{split}
\mathfrak{N}_x(v) & = \int_{\M} \Pi_{\exp_{x_v}}^x \log_{x_v}(y)\: \mu(dy)
= \int_{\M} l_x(v, \log_x(y) )\: \mu(dy) \\
& = \int_\M  (\log_x(y)-v)\: \mu(dy) 
+ \frac{1}{6} \int_\M R(\log_x(y) ,v)(v-2 \log_x(y))\: \mu(dy) \\
&\quad + \int_\M  \frac{1}{24}(\nabla_v R)(\log_x(y),v) (2v-3\log_x(y)) \: \mu(dy) \\
&\quad + \int_\M  \frac{1}{24}(\nabla_{\log_x(y)} R)(\log_x(y),v) (v-2\log_x(y))\: \mu(dy) + \int_\M O(5)\: \mu(dy) .
\end{split}
\]
In the integral of the first terms of the series, we recognize the first moment and the normalization (zeroth moment) of the distribution:
$\int_\M (\log_x(y)-v)\: \mu(dy) = \MM_1 - v.$
For the second order term, we take the curvature tensor out of the integral thanks to its mutlilinearity, and we see the first and second moment appearing:
\[
\int_\M  R(\log_x(y) ,v)(v-2 \log_x(y))\: \mu(dy) 
%= \int_\M  R(\log_x(y) ,v) v \mu(dy) -2 \int_\M  R(\log_x(y) ,v)\log_x(y) \mu(dy) 
= R(\MM_1,v) v -2 R(\contractvar ,v)\contractvar \contract \MM_2.
\]
In this formula, the notation $R(\contractvar, v)\contractvar \contract \MM_2$ denotes the contraction of the $(2,0)$ symmetric moment tensor with the $(1,3)$ curvature tensor $R$ along the axes specified by the bullets. In coordinates, this writes 
$ [ R(\contractvar ,v)\contractvar \contract \MM_2 ]^a = R^a_{bcd}\: v^d \: \MM_2^{bc}$.  In further computation, we will also use contractions like $
\nabla_\contractvar R(\contractvar ,v)\contractvar \contract \MM_3$ or 
$\MM_2 \contractcirc \nabla_{\contractvar} R(\contractvar, \contractvarcirc) \contractvarcirc \contract \MM_2  $. 
For the 4th order terms, we get:
\[
\int_\M  (\nabla_v R)(\log_x(y),v) (2v-3\log_x(y))\: \mu(dy) 
%= 2 \int_\M  (\nabla_v R)(\log_x(y),v) v \mu(dy) -3 \int_\M  (\nabla_v R)(\log_x(y),v)\log_x(y) \mu(dy)
= 2 (\nabla_v R)(\MM_1,v) v - 3 (\nabla_v R)(\contractvar,v)\contractvar \contract \MM_2 ,
\]
\[
\int_\M (\nabla_{\log_x(y)} R)(\log_x(y),v) (v-2\log_x(y))\: \mu(dy) 
%=  \int_\M (\nabla_{\log_x(y)} R)(\log_x(y),v) v \mu(dy) -2 \int_\M (\nabla_{\log_x(y)} R)(\log_x(y),v) \log_x(y) \mu(dy) = 
	=	(\nabla_{\contractvar} R)(\contractvar,v)v \contract\MM_2 -2 (\nabla_{\contractvar} R)(\contractvar,v) \contractvar \contract \MM_3 .
\]
Last but not least, $O(5)$ is a homogeneous polynomial of order at least 5 in $v$ and $\log_x(y)$. Both vectors have a norm less than the diameter $\varepsilon$  of the distribution, so that $\int_\M O(5) \: \mu(dy) = O(\varepsilon^5)$. We finally obtain the following expansion.  

\begin{theorem}[The recentered tangent mean map expansion]
\label{thm:RecenteredMeanMapExpansion}
Let $\mu$ be a probability distribution with support in a convex open subset $\U$ of diameter $\varepsilon$ of an affine manifold $\M$ satisfying the KKC or ALC conditions. 
The expansion of the recentered tangent mean map around a point $x\in \U$ with respect to the vector $v\in V_x = \log_x( \U ) \subset T_x\M$ is:
\begin{equation}
\begin{split}
\mathfrak{N}_x(v) & =
\MM_1 - v 
+ \frac{1}{6} R(\MM_1,v) v - \frac{1}{3} R(\contractvar ,v)\contractvar \contract \MM_2 \\
& \quad + \frac{1}{12} (\nabla_v R)(\MM_1,v) v -  \frac{1}{8} (\nabla_v R)(\contractvar,v)\contractvar \contract \MM_2  \\
& \quad + \frac{1}{24} (\nabla_{\contractvar} R)(\contractvar,v) v \contract \MM_2 -\frac{1}{12} (\nabla_{\contractvar} R)(\contractvar,v) \contractvar \contract \MM_3 
+ O(\epsilon^5)
\label{eq:RecenteredMeanMapExpansion}
\end{split}
\end{equation}
\end{theorem}

\subsection{Solving for the zero of the recentered mean map}

Now that we have the expansion of the recentered mean map, the last step is to solve for the value of $v$ that zeros it out. The second order approximation of $\mathfrak{N}_x(v)=0$ is obviously $v= \MM_1 + O(\varepsilon^3)$ since the curvature terms are of order three. Now, to identify the third order term $v_3$, we plug the value $\MM_1 +v_3 +O(\varepsilon^4)$ in Eq.(\ref{eq:RecenteredMeanMapExpansion}):
\[
\begin{split}
\mathfrak{N}_x(\MM_1 +v_3 +O(\varepsilon^4)) & = - v_3 
+ \frac{1}{6} R(\MM_1,\MM_1) \MM_1 
- \frac{1}{3} R(\contractvar ,\MM_1)\contractvar \contract \MM_2 
+O(\varepsilon^4).
\end{split}
\]
The first term cancels out thanks to the skew-symmetry of the curvature tensor, so that 
we find that $v =\MM_1 - \frac{1}{3} R(\contractvar ,\MM_1)\contractvar \contract \MM_2 +O(\varepsilon^4)$ is required to cancel out the recentered mean map. The fourth order give:
\[
\begin{split}
&\mathfrak{N}_x\left(\MM_1 - \frac{1}{3} R(\contractvar ,\MM_1)\contractvar \contract \MM_2 +v_4 +O(\varepsilon^5) \right) =  - v_4
%% this term is of order 5 : &\qquad   + \frac{1}{6} R(\MM_1, - \frac{1}{3} R(\contractvar ,\MM_1)\contractvar : \MM_2 +v_4) (\MM_1 - \frac{1}{3} R(\contractvar ,\MM_1)\contractvar \contract \MM_2 +v_4) \\
  + \frac{1}{12} (\nabla_{\MM_1} R)(\MM_1,\MM_1) \MM_1 
\\ & \quad -  \frac{1}{8} (\nabla_{\MM_1} R)(\contractvar,\MM_1)\contractvar \contract \MM_2  
 + \frac{1}{24} (\nabla_{\contractvar} R)(\contractvar,\MM_1) \MM_1 \contract \MM_2 -\frac{1}{12} (\nabla_{\contractvar} R)(\contractvar, \MM_1) \contractvar \contract \MM_3 
+ O(\epsilon^5) .
\end{split}
\]
Removing the vanishing term due to  the skew-symmetry of the curvature tensor, we find:
\[
\begin{split}
 v_4 =   -  \frac{1}{8} (\nabla_{\MM_1} R)(\contractvar,\MM_1)\contractvar \contract \MM_2  
  + \frac{1}{24} (\nabla_{\contractvar} R)(\contractvar,\MM_1) \MM_1 \contract \MM_2 -\frac{1}{12} (\nabla_{\contractvar} R)(\contractvar, \MM_1) \contractvar \contract \MM_3 
+ O(\epsilon^5).
\end{split}
\]

\begin{theorem}
\label{thm:FrechetMeanField}
Let $\mu$ be a probability distribution with support in a convex open subset $\U$ of diameter $\varepsilon$ of an affine manifold $\M$ satisfying the KKC or ALC conditions. 
The vector field that points from a point of $\U$ to the mean $\bar x$ characterized by $\MM_1(\bar x) =0$ has the following Taylor expansion:
\begin{equation}
\label{eq:LogFrechet}
\begin{split}
\log_x(\bar x) = &\: 
\MM_1  - \frac{1}{3} R(\contractvar, \MM_1)\contractvar \contract \MM_2 
+ \frac{1}{24} \nabla_\contractvar R(\contractvar,\MM_1)\MM_1 \contract \MM_2 
\\ & -\frac{1}{8} \nabla_{\MM_1} R(\contractvar, \MM_1) \contractvar \contract \MM_2
-\frac{1}{12} \nabla_\contractvar R(\contractvar, \MM_1) \contractvar \contract\MM_3 + O(\epsilon^5).
\end{split}
\end{equation} 
\end{theorem}

In coordinates, 
we have $[R(u,v)w]^a = R^a_{bcd} w^b u^c v^d$  
% $[R(\contractvar, \MM_1)\contractvar \contract \MM_2]^a = R^a_{bcd} \MM_2^{bc} \MM_1^d$ 
and $[\nabla_u R(v,w)z]^a = \nabla_b R^a_{cde} z^c v^d w^e = R^a_{cde;b} z^c v^d w^e$, 
so that this formula writes:
\[
\begin{split}
\log_x(\bar x)^a = &\:
\MM_1^a  -\frac{1}{3} % R(\contractvar, \MM_1)\contractvar \contract \MM_2
                        R^a_{bcd} \MM_2^{bc} \MM_1^d 
+ \frac{1}{24} %\nabla_\contractvar R(\contractvar,\MM_1)\MM_1 \contract \MM_2 
                \nabla_e R^a_{bcd} \MM_2^{ce} \MM_1^b \MM_1^d 
\\ & -\frac{1}{8} %\nabla_{\MM_1} R(\contractvar, \MM_1) \contractvar \contract \MM_2
                  \nabla_e R^a_{bcd} \MM_1^e \MM_2^{bc} \MM_1^d 
-\frac{1}{12} % \nabla_\contractvar R(\contractvar, \MM_1) \contractvar \contract\MM_3 
              \nabla_e R^a_{bcd} \MM_1^d  \MM_3^{bce}
+ O(\epsilon^5),
\end{split}
\]
The order 3 of this expression was earlier derived by \cite{darling_geometrically_2000}[Theorem 3.2] with an opposite sign convention for the curvature (see appendix \ref{sec:Darling} for the details of the notations equivalences).
In this work, we add the order 4, which will turn out to be crucial for establishing the bias on the empirical Fr\'echet mean, and the method to get higher orders.

\section{Non-asymptotic high concentration expansion of the empirical mean}
\label{sec:IIDSample}

Let $\mu$ be a distribution on $\M$ satisfying the KKC or ALC conditions, with mean $\bar x$. 
An IID $n$-sample $X_n = \{x_1, \ldots x_n\} \in \M^n$  may be identified with the  empirical distribution $X_n \simeq \frac{1}{n} \sum_i \delta_{x_i}$. To study the law of its empirical mean $\bar x_n$,  we focus in this paper on its expected moments.
 To simplify the notations, we denote in this section $\MM_k = \MM_k\mu$ the moments of the underlying distribution $\mu$ and $\XX_k^n = \MM_k X_n$ the $k$-th moments of the $n$-sample $X_n$. 

Since $X_n$ is a distribution, we can use directly Theorem \ref{thm:FrechetMeanField} to find that the location of the empirical Fr\'echet mean $\bar x_n$ with respect to a point $x$:
\begin{equation}
\label{eq:LocEmpiricalFrechetMean}
\begin{split}
\log_x(\bar x_n) = \: & 
\XX^n_1  - \frac{1}{3} R(\contractvar, \XX^n_1)\contractvar \contract \XX^n_2 
+ \frac{1}{24} \nabla_\contractvar R(\contractvar,\XX^n_1)\XX^n_1 \contract \XX^n_2 
\\ & -\frac{1}{8} \nabla_{\XX^n_1} R(\contractvar, \XX^n_1) \contractvar \contract \XX^n_2
-\frac{1}{12} \nabla_\contractvar R(\contractvar, \XX^n_1) \contractvar \contract\XX^n_3 + O(\epsilon^5) .
\end{split}
\end{equation}

This formula makes use of the empirical moment 
$\XX_k^n(x)    % = \MM_k X_n(x) 
= \frac{1}{n} \sum_{i=1}^n \lcp{x x_i} \otimes \ldots \otimes \lcp{x x_i}$ of the $n$-sample $X_n$. 
 The expectation of this tensor when the $x_i$'s are IID with law $\mu$ is simply 
\[
\expect{\XX_k^n(x) } =  \frac{1}{n} \sum_{i=1}^n \expect{\lcp{x x_i} \otimes \ldots \otimes \lcp{x x_i}} 
=  \frac{1}{n} \sum_{i=1}^n \int_{\M} \lcp{x x_i} \otimes \ldots \otimes \lcp{x x_i}
\: \mu(dx_i) = \MM_k (x). 
\]
However, the expectation of the tensor product of moments is more complex than the simple contraction of their expectation. The reason is that the random variables $x_i$ and $x_j$ appearing in two empirical moments are independent if $i \not = j$ but not if $i=j$, in which case a higher order moments appears.

\subsection{Expectation of tensor product of empirical moments}
\label{sec:ExpectEmpiricalMoments}
In this section, we keep the sums over the variables $i,j,k$ explicit to stress that these are the indices of the data points and not covariant indices of vectors.
For the product of the first and second moments appearing above, we have: 
\[
\begin{split}
 \expect{  (\XX_1^n) \otimes (\XX_2^n) }  
  & =   \frac{1}{n^2} \sum_{i,j} \expect{   \lcp{x x_i} \otimes  \lcp{x x_j} \otimes \lcp{x x_j} }\\
	& =   \frac{1}{n^2} \sum_{i, j\not = i} \expect{ \lcp{x x_i} \otimes  \lcp{x x_j} \otimes \lcp{x x_j} }
	 +   \frac{1}{n^2} \sum_{i} \expect{ \lcp{x x_i} \otimes  \lcp{x x_i} \otimes \lcp{x x_i} } \\
 & = \frac{n-1}{n} \MM_1 \otimes \MM_2  + \frac{1}{n}  \MM_3 
\end{split}
\]
This computation can be generalized to the expectation of the product of two moments of any order: 
\[
\begin{split}
 \expect{  \XX_p^n \otimes \XX_q^n}  
	& = \frac{n-1}{n} \MM_p \otimes \MM_q+ \frac{1}{n}  \MM_{p+q}
\end{split}
\]

For three moments, we have:
%\[
%\begin{split}
 %\sum_{i, j, k } %\expect
%{ ( \otimes_p v_i) ( \otimes_q v_j)  (\otimes_r v_k) } =
	%&   \sum_{i, j \not = i, k \not = i,j } %\expect
	%{ ( \otimes_p v_i) ( \otimes_q v_j)  (\otimes_r v_k) } 
	     	%+ \sum_{i, j \not = i, k = i } %\expect
				%{ ( \otimes_{(p+r)} v_i) ( \otimes_q v_j)  } \\
	%&   	+  \sum_{i, j \not = i, k = j } %\expect
	%{ ( \otimes_{p} v_i) ( \otimes_{(q+r)} v_j)  }
	    	%+  \sum_{i, j = i, k \not = i} %\expect
				%{ ( \otimes_{(p+q)} v_i) (\otimes_r v_k) } \\
	%&   	+  \sum_{i = j = k } %\expect
	%{ ( \otimes_{(p+q+r)} v_i) },
%\end{split}
%\]
%so that:
\[
\begin{split}
 \expect{  \XX_p^n \otimes \XX_q^n \otimes \XX_r^n }  
	 = & \frac{1}{n^3} \sum_{i, j, k } \expect{ ( \lcp{x x_i}^{\otimes p}) ( \lcp{x x_j}^{\otimes q } )  (\lcp{x x_k}^{\otimes r }) }\\
	 = & \frac{(n-1)(n-2)}{n^2} \MM_p  \otimes \MM_q \otimes \MM_r + \frac{1}{n^2} \MM_{p+q+r} \\
	   &  + \frac{(n-1)}{n^2} ( \MM_{p+q} \otimes \MM_r + 
		%\MM_{p+r} \otimes \MM_q 
		(\contractvar \otimes \MM_q \otimes \contractvar) \contract \MM_{p+r}
		+ \MM_p \otimes \MM_{q+r} ).
\end{split}
\]

\subsection{First moment of the empirical  mean}

Now we can come back to the location of the empirical Fr\'echet mean $\bar x_n$. To obtain its first moment, we take the expectation of $\log_x(\bar x_n)$ with respect to the joint law $\mu^{\otimes n}$ of the sample. Using Eq.\ref{eq:LocEmpiricalFrechetMean}, this gives:
%theorem \ref{thm:FrecherMeanField}: 
\[
\begin{split}
 \expect{ \log_x(\bar x_n) } = &
\expect{ \XX^n_1}  - \frac{1}{3}  \expect{ R(\contractvar, \XX^n_1)\contractvar \contract \XX^n_2 }
+ \frac{1}{24}  \expect{ \nabla_\contractvar R(\contractvar,\XX^n_1)\XX^n_1 \contract \XX^n_2 }
\\ & -\frac{1}{8}  \expect{ \nabla_{\XX^n_1} R(\contractvar, \XX^n_1) \contractvar \contract \XX^n_2 }
-\frac{1}{12}  \expect{ \nabla_\contractvar R(\contractvar, \XX^n_1) \contractvar \contract\XX^n_3} + O(\epsilon^5) .
\end{split}
\]
Then we expand the expectation of the tensor product of empirical moments using the rules of Section \ref{sec:ExpectEmpiricalMoments}:  
\[
\begin{split}
 \expect{ \log_x(\bar x_n) } = &
\MM_1  - \frac{n-1}{3n}  R(\contractvar, \MM_1)\contractvar \contract \MM_2 
- \frac{1}{3n}  R(\contractvar, \contractvar)\contractvar \contract \MM_3 
\\ & + \frac{1}{24 n^2}  \nabla_\contractvar R(\contractvar,\contractvar)\contractvar \contract \MM_4 
+ \frac{(n-1)(n-2)}{24 n^2}  \nabla_\contractvar R(\contractvar,\MM_1)\MM_1 \contract \MM_2 
\\ & + \frac{(n-1)}{24 n^2} \left( 
\nabla_\contractvar R(\contractvar,\contractvarcirc)\contractvarcirc \contract \MM_2 \contractcirc \MM_2
+ \nabla_\contractvar R(\contractvar,\MM_1)\contractvar \contract \MM_3 
+ \nabla_\contractvar R(\contractvar,\contractvar)\MM_1 \contract \MM_3  \right)
\\ & -\frac{1}{8 n^2}  \nabla_{\contractvar} R(\contractvar, \contractvar) \contractvar \contract \MM_4 
 -\frac{(n-1)(n-2)}{8 n^2}  \nabla_{\MM_1} R(\contractvar, \MM_1) \contractvar \contract \MM_2 
\\ & -\frac{(n-1)}{8 n^2}  \left( 
	    \nabla_{\contractvarcirc} R(\contractvar, \contractvarcirc) \contractvar \contract \MM_2 \contractcirc \MM_2
	+ \nabla_{\contractvar} R(\contractvar, \MM_1) \contractvar \contract \MM_3
	+ \nabla_{\MM_1} R(\contractvar, \contractvar) \contractvar \contract \MM_3
	\right)
\\ & -\frac{(n-1)}{12 n} 	\nabla_\contractvar R(\contractvar, \MM_1) \contractvar \contract\MM_3
	-\frac{1}{12 n} \nabla_\contractvar R(\contractvar, \contractvar) \contractvar \contract\MM_4
 + O(\epsilon^5) .
\end{split}
\]
We now remove the terms that vanish because they are the contraction of a symmetric tensor $\MM_k$ with skew symmetric indices of the curvature tensor and we factor similar terms (sometimes changing sign thanks to the skew-symmetry again). We finally find:

\begin{theorem}[First moment of the empirical mean] 
\label{thm:EmpiricalFrechetMeanField}
Let $\mu$ be a probability distribution satisfying the KKC or ALC conditions with support of diameter less than $\varepsilon$ and mean $\bar x$. 
The first moment of the empirical Fr\'echet mean $\bar x_n$ of an IID $n$-sample is:
\begin{equation*}
\label{eq:EmpiricalFrechetMeanField}
\begin{split}
 \expect{ \log_x(\bar x_n) } = &\: \MM_1  - \frac{n-1}{3n}  R(\contractvar, \MM_1)\contractvar \contract \MM_2 
\\& + \frac{(n-1)(n-2)}{24 n^2} \left( \nabla_\contractvar R(\contractvar,\MM_1)\MM_1 \contract \MM_2 
	-3 \nabla_{\MM_1} R(\contractvar, \MM_1) \contractvar \contract \MM_2 \right)
\\ &    +\frac{(n-1)}{12 n^2} \left(  2 \nabla_{\contractvarcirc} R(\contractvarcirc, \contractvar) \contractvar \contract \MM_2 \contractcirc \MM_2 - \left(1+n\right)\nabla_\contractvar R(\contractvar,\MM_1)\contractvar \contract \MM_3 \right)
 + O(\epsilon^5) .
\end{split}
\end{equation*}
\end{theorem}

At the mean $\bar x$ of the distribution $\mu$, the first moment $\MM_1 (\bar x)=0$  vanishes, so that we end up with the estimation of how the empirical mean $\bar x_n$ deviates from the mean $\bar x$ of the underlying distribution $\mu$ in expectation (i.e. the bias):
\begin{equation}
\label{eq:EmpiricalFrechetMeanBias}
 \expect{ \log_{\bar x}(\bar x_n) } =  
   \frac{(n-1)}{6 n^2}  \nabla_{\contractvar} R(\contractvar, \contractvarcirc) \contractvarcirc \contract \MM_2 \contractcirc \MM_2 
 + O(\epsilon^5) .
\end{equation}

\subsection{Second moment of the empirical  mean}

Let $\Sigma_x^n = \log_x(\bar x_n) \otimes \log_x( \bar x_n)$ be the tensor product of the empirical mean with itself. From Eq. \ref{eq:LocEmpiricalFrechetMean}, we can write it as:
\[
\begin{split}
%\log_x(\bar x_n) \otimes \log_x( \bar x_n) = &
\Sigma_x^n = &
\: \XX^n_1 \otimes \XX^n_1 - \frac{1}{3} (R(\contractvar, \XX^n_1)\contractvar \contract \XX^n_2) \otimes \XX^n_1 - \frac{1}{3} \XX^n_1 \otimes  (R(\contractvar, \XX^n_1)\contractvar \contract \XX^n_2 )  + O(\epsilon^5).
\end{split}
\]
Taking the expectation to obtaining the covariance field, we have to expand the expression of the tensor product of empirical moments:
\[
\expect{ \XX^n_1 \otimes \XX^n_1 } = \frac{1}{n} \MM_2 + \frac{n-1}{n}\MM_1 \otimes \MM_1 ;
\]
\[
\begin{split}
\expect{ \XX^n_1 \otimes  (R(\contractvar, \XX^n_1)\contractvar \contract \XX^n_2 )} = &
 \: \frac{n-1}{n^2} \Big( 
				\MM_2 \contractcirc \contractvarcirc \otimes R(\contractvar, \contractvarcirc)\contractvar \contract \MM_2 
				+ \contractvar \otimes  R(\contractvar, \MM_1)\contractvar \contract \MM_3
 \Big)
 \\ & +\frac{(n-1)(n-2)}{n^2} \MM_1 \otimes  R(\contractvar, \MM_1)\contractvar \contract \MM_2.
\end{split}
\]

Recombining the terms, we get:
\[
\begin{split}
\expect{ \Sigma_x^n } = &
\MM_2 + \frac{n-1}{n}\MM_1 \otimes \MM_1 
-\frac{(n-1)}{3 n^2} 
				\MM_2 \contractcirc ( \contractvarcirc \otimes R(\contractvar, \contractvarcirc)\contractvar + R(\contractvar, \contractvarcirc)\contractvar \otimes \contractvarcirc  )\contract \MM_2 
\\ &   - \frac{(n-1)(n-2)}{3 n^2} \Big(
	R(\contractvar, \MM_1)\contractvar \contract \MM_2 \otimes \MM_1 
	+ \MM_1 \otimes  R(\contractvar, \MM_1)\contractvar \contract \MM_2 \Big)
\\& - \frac{n-1}{3 n^2} \Big( 
          R(\contractvar, \MM_1)\contractvar \contract \MM_3 \otimes \contractvar
				+ \contractvar \otimes  R(\contractvar, \MM_1)\contractvar \contract \MM_3
 \Big)     + O(\epsilon^5)
\end{split}
\]
At the mean $\bar x$ of the distribution $\mu$, the first moment $\MM_1 (\bar x)=0$  vanishes, so that we end up with the following estimation of the covariance matrix of  the empirical mean:
\begin{equation}
\label{eq:EmpiricalFrechetMeanCov}
 \expect{ \log_{\bar x}(\bar x_n) \otimes \log_{\bar x}(\bar x_n) } =  
   \frac{1}{n}\MM_2  -\frac{(n-1)}{3 n^2} 
				\MM_2 \contractcirc ( \contractvarcirc \otimes R(\contractvar, \contractvarcirc)\contractvar + R(\contractvar, \contractvarcirc)\contractvar \otimes \contractvarcirc  )\contract \MM_2 
 + O(\epsilon^5) .
\end{equation}
In coordinates, this write perhaps more simply:
\[
 \expect{ \log_{\bar x}(\bar x_n) \otimes \log_{\bar x}(\bar x_n) }^{ab} =  
   \frac{1}{n}\MM_2^{ab}  -\frac{(n-1)}{3 n^2} 
\MM_2^{cd} ( \MM_2^{ae}   R^b_{cde} + R^a_{cde} \MM_2^{be})
 + O(\epsilon^5) .
\]

\begin{theorem}[First moments of the empirical mean] 
%$ $\\
\label{thm:MomentsEmpiricalMean}
Let $\mu$ be a probability distribution satisfying the KKC or ALC conditions with support of diameter less than $\varepsilon$ and mean $\bar x$. 
We denote $\MM_k = \MM_k\mu (\bar x)$ its $k$-order moment at the mean $\bar x$. By definition, we have $ \MM_0=1$ and  $ \MM_1 =0$. 
Let $X_n = \{x_1, \ldots x_n\} \in \M^n$ be an IID $n$-sample of this distribution. The empirical mean $\bar x_n$ of this sample is also unique. Its expected first moment  at $\bar x$, the bias $\text{Bias}(\bar x_n) = \expect{ \log_{\bar x}(\bar x_n) }$, is
\begin{equation}
\label{eq:EmpiricalFrechetMeanBiasThm}
 \text{Bias}(\bar x_n) = %\expect{ \log_{\bar x}(\bar x_n) } =  
  \textstyle  \frac{1}{6 n}\left(1-\frac{1}{n}\right)  \MM_2 \contractcirc \nabla_{\contractvar} R(\contractvar, \contractvarcirc) \contractvarcirc \contract \MM_2  
 + O\left(\epsilon^5 \right) . 
\end{equation}
Its expected second moment  at $\bar x$,
the covariance $\text{Cov}(\bar x_n) = 
 \expect{ \log_{\bar x}(\bar x_n) \otimes \log_{\bar x}(\bar x_n) }$, is:
\begin{equation}
\label{eq:EmpiricalFrechetMeanCovThm}
\text{Cov}(\bar x_n) = 
 % \expect{ \log_{\bar x}(\bar x_n) \otimes \log_{\bar x}(\bar x_n) } =  
  \textstyle \frac{1}{n} \left( \MM_2  -\frac{1}{3}\left(1-\frac{1}{n}\right)
				\MM_2 \contractcirc ( \contractvarcirc \otimes R(\contractvar, \contractvarcirc)\contractvar + R(\contractvar, \contractvarcirc)\contractvar \otimes \contractvarcirc  )\contract \MM_2   \right)
 + O\left(\epsilon^5\right) .  
\end{equation}
In coordinates, this writes
\begin{eqnarray*} 
\text{Bias}(\bar x_n)^a & = &  
  \textstyle \frac{1}{6 n}\left(1-\frac{1}{n}\right)  \nabla_b R^a_{cde} \MM_2^{ce} \MM_2^{bd} 
 + O\big(\epsilon^5) \\ 
\text{Cov}(\bar x_n)^{ab} & = & \textstyle \frac{1}{n}\big(  \MM_2^{ab} - \frac{1}{3}\left(1-\frac{1}{n}\right)
\MM_2^{cd} ( \MM_2^{ae}   R^b_{cde} + R^a_{cde} \MM_2^{be}) \big) + O\big(\epsilon^5) . 
\end{eqnarray*}
\end{theorem}

Thus, in general manifolds, there is a bias on the empirical mean of order $1/n$ in the number of samples and of order 4 in $\varepsilon$. Asymptotically, this bias disappear with the number of samples. However, keeping a fixed number of samples,
% and a fixed covariance matrix, 
this bias increases and could possibly blow up (outside the domain of validity of our assumptions) when we approach a singularity with an unbounded curvature, in which case the gradient of the curvature has to become large. Of course, the KKC or ALC conditions 
on the support of the distribution only holds at a sufficient distance of such a singularity for a fixed covariance, which limits the conclusion that we may draw from this trend.
Interestingly, such a bias is not visible in symmetric spaces since the curvature tensor is covariantly constant in these manifolds.

The covariance matrix of the empirical mean also has a curvature correction term modulating (accelerating or decelerating) the convergence. In order to better understand the effect of curvature, we study in the next section how this relates to the Bhattacharya-Patrangenaru central limit theorem in Riemannian manifolds, and in Section \ref{sec:ModulartionConstantCurvature} how the formulas simplify for constant curvature spaces.

\section{Asymptotic covariance of the empirical Riemannian mean}
\label{sec:BP-CLT}

The Bhattacharya-Patrangenaru central limit theorem (BP-CLT) for sample Fr\'echet means in manifolds \cite[Theorem 2.1 p. 1230]{bhattacharya_large_2005} is a general CLT valid for non-Riemannian twice differentiable distances.
%, under quite technical assumptions. 
We consider here the case of the intrinsic Fr\'echet mean on a Riemannian manifold \cite[Theorem 2.2 p. 1231]{bhattacharya_large_2005}.

\subsection{The Bhattacharya-Patrangenaru CLT}
We first align our notations to the ones of \cite{bhattacharya_large_2005}. Recall that the support of $\mu$ is included in $\U \subset B(x,r)$ with $r < r^*$ by the KKC condition.  Because we use a normal coordinate system, a point $y \in \U$  is parametrized by $v=\log_x(y) \in V_x = \log_x(\U)$, so that our $v$ is the $\theta$ for BP and the BP function  $\phi(y)$ is $ \log_x(y)$ for us.  With our notation $x_v = \phi\inv(v) = \exp_x(v)$, the pullback of the Riemannian distance to the chart writes: $\left( \rho^{\phi}(u,v) \right)^2  = \dist({x_u}, {x_v})^2 = \| l_x(u,v)\|^2_x$.   The BP function $\Psi(u ; v) = \left( \frac{\partial}{\partial v} \dist( x_u , {x_v})^2 \right)^\top$ is the Euclidean gradient of this pulled back squared distance with respect to $v$. 
The BP matrix $D_v \Psi(u;v) $ is second order derivative of the pulled back squared distance $v \rightarrow \dist({x_u}, {x_v})^2$ in our chart. In the original BP-CLT theorem, the covariance of the normal law is $\Lambda\inv \: \Sigma \: \Lambda\inv$, where $\Sigma$ is the covariance of $\Psi(u ; v)$ and the matrix $\Lambda = \expect{D_v \Psi(\log_{x}(x_i) ; v)}$ is the  expected value of the second order derivative under the law $\mu$ of the sample $x_i$. 

In a normal coordinate system centered at $\bar x$, the Christoffel symbols and their derivatives vanish at $\bar x$ so that the standard differential corresponds to the Riemannian gradient and the standard second order derivative corresponds to the Riemannian Hessian. Let $d_y^2(x) = \dist^2(y,x)$. For $y \not \in \text{Cut}(x)$, we have $\nabla d_y^2(x) = -2 \log_{x}(y)$ and the Hessian is $H_x(y) = \nabla^2 d_y^2(x) = -2 D_x \log_{x}(y)$ (see for instance \cite[appendix A]{pennec2018}). Thus, we have $\Psi(u ; 0) = \nabla d_{x_v}^2(x_u)|_{v=0} = -2 \log_{\bar x}(x_u) = -2 u$ and its covariance is simply $\Sigma = 4 \int_\M \log_{\bar x}(y) \log_{\bar x}(y)^\top \: \mu(dy) = 4 \MM_2$. For the second order derivatives, we have $D_v \Psi(\log_{\bar x}(y) ; v)|_{v=0} = H_{\bar x}(y)$ and its expected value for $y$ following the law $\mu$ is $\Lambda = \bar H = \int_\M H_{\bar x}(y) \: \mu(dy)$. 
These alignments of notations lead to the following formulation of the BP-CLT for intrinsic means. 

%rewritten with our notations in a normal coordinate system:
\begin{theorem}[BP-CLT for intrinsic sample means \cite{bhattacharya_large_2005}] 
\label{BP-CLT}
Let $\mu$ is a probability distribution on a Riemannian manifold $(\M,g)$ of Fr\'echet mean $\bar x$ and covariance matrix $\MM_2 = \MM_2 \mu(\bar x)$ whose support is included in a regular geodesic ball $B(x_0,r)$ satisfying the KKC conditions. % (definition \ref{KKCC}). 

Then (a) the empirical Fr\'echet mean $\bar x_n$ of a sample $X_n = \frac{1}{n} \sum_{i=1}^n \delta_{x_i}$ is a consistent estimator of the population Fr\'echet mean $\bar x$ of $\mu$ and (b) the random variable $\sqrt{n} \log_{\bar x} (\bar x_n) \in T_{\bar x} \M$ converges in law to  a normal distribution of mean 0 and covariance $4 \bar H\inv \: \MM_2 \: \bar H\inv$, where the matrix $\bar H$ is the expectation of the Riemannian Hessian of the squared distance $\dist( . , y)^2$ according to the distribution $\mu$.
\end{theorem}

In order to compare our non-asymptotic high concentration expansion (Theorem \ref{thm:MomentsEmpiricalMean}) to the BP-CLT, we first observe that the consistency of the empirical Fr\'echet mean of the BP-CLT states that the limit of the bias in Eq.(\ref{eq:EmpiricalFrechetMeanBiasThm}) should vanish when $n$ goes to infinity. This is indeed the case since 
 $\text{Bias}(\bar x_n) = \expect{\log_{\bar x}(\bar x_n)}$ is of order $1/n$. 
To compare the  covariance, we essentially need to compute the Hessian of the squared Riemannian distance $H_x(y)$, integrate it to get its expectation $\bar H$, and show that $4 \bar H \inv \MM_2 \bar H\inv$ has the same Taylor expansion as in  Eq.(\ref{eq:EmpiricalFrechetMeanCovThm}).

\subsection{Hessian of the squared Riemannian distance}
The Taylor expansion of $H_x(x_u)$ has been established in Eq.(1) of \cite{pennec2018} using Brewin Taylor expansion\footnote{We indicate here the order $O(4)$ corresponding to the notations of the current paper since the third order of Brewin expansions has a different meaning which includes a conformal factor accounting for curvature.}:
\[
\frac{1}{2} [H_x(x_u)]^a_b = - [D_x \log_x(x_u)]^a_b = \delta^a_b +\frac{1}{3} R^a_{cdb} u^c u^d +\frac{1}{12} \nabla_c R^a_{deb} u^c u^d u^e + O(4).
\]
We verify below that we obtain the same formulation using the Gavrilov's Taylor expansions used in the current paper. 
The Taylor expansion of the squared distance in a normal coordinate system at $x$ is given by Eq.(\ref{eq:TaylorSquareDist}): \[
\dist^2(x_v, x_u) =  (u - v)^a   (u - v)_a  + \frac{1}{3} R_{abcd} v^a u^b u^c v^d +\frac{1}{12} \nabla_e R_{abcd} v^a u^b u^c v^d (u+v)^e  +O(6).
\]
  Taking the derivative with respect to the coordinate $\alpha$ of $v$ and raising the index $\alpha$, we get the gradient:
\[
%\partial_{v_\alpha} \dist^2(x_v, x_u)  
\Psi(u,v) = -2 (u-v)^\alpha  
+ \frac{2}{3} R^{\alpha}_{bcd} u^b u^c v^d 
+ \frac{1}{6} \nabla_e R^{\alpha}_{bcd} u^b u^c v^d (u+v)^e  
+ \frac{1}{12} \nabla^{\alpha} R_{ebcd} v^e u^b u^c v^d  
+O(5).
\]
Taking the derivative with respect to the coordinate $\beta$ of $v$, we have:
\[
 D_\beta \Psi(u,v) = 2 \delta^{\alpha}_{\beta}   
+ \frac{2}{3} R^{\alpha}_{bc\beta} u^b u^c 
+ \frac{1}{6} \nabla_e R^{\alpha}_{bc\beta} u^b u^c (u+v)^e  
+ \frac{1}{6} \nabla_{\beta} R^{\alpha}_{ecd} u^e u^c v^d    
+ \frac{1}{6} \nabla^{\alpha} R_{ebc\beta} v^e u^b u^c  
 + O(4)
\]
The value at $v=0$ is $\frac{1}{2}[H_x(x_u)]^{\alpha}_{\beta} =  \delta^{\alpha}_{\beta}   + \frac{1}{3} R^{\alpha}_{bc\beta} u^b u^c 
 +\frac{1}{12} \nabla_e R^{\alpha}_{bc\beta} u^b u^c u^e  + O(4),$
in accordance with the previous expansion of \cite{pennec2018}. 
The expectation of this Hessian matrix when $x_u$ has probability $\mu$ is:
 \[
\frac{1}{2} \bar H^{a}_{b} =  \delta^a_b   + \frac{1}{3} R^a_{dcb} \MM_2^{dc}
 +\frac{1}{12} \nabla_e R^a_{dcb} \MM_3^{dce}  + O(\varepsilon^4)
\]

\subsection{Covariance matrix of the BP-CLT}

The inverse of the matrix $\bar H$ can be determined by identification in $[\bar H\inv]^{a}_{b} [\bar H]^b_c = \delta^a_c$: 
 \[
 [\bar H\inv]^{a}_{b} =  \frac{1}{2} \left( \delta^a_b   - \frac{1}{3} R^a_{dcb} \MM_2^{dc}
 - \frac{1}{12} \nabla_e R^a_{dcb} \MM_3^{dce} \right)  + O(\varepsilon^4).
\]
Finally, we get
\[
\begin{split}
4 [\bar H \inv \MM_2 \bar H\inv]^{ab} & = 
4 [\bar H \inv]^a_c [\MM_2]^{cd} [\bar H\inv]^b_d 
	= \MM_2^{ab} - \frac{1}{3} \MM_2^{ef}  \left( R^a_{efc} \MM_2^{cb} + \MM_2^{ad} R^b_{efd} \right)
	+ O(\varepsilon^5),
	\end{split}
\] 
which is in accordance with the term in $1/n$ of Eq.\ref{eq:EmpiricalFrechetMeanCovThm} of Theorem \ref{thm:MomentsEmpiricalMean}.

\section{Modulation of the rate of convergence in space forms}
\label{sec:ModulartionConstantCurvature}

In order to verify experimentally our predicted bias and rate of convergence, it is better to have explicit formulas for geodesics, which are known only for symmetric spaces. In these spaces, the curvature is covariantly constant, so that the empirical mean $\bar x_n$ of an IID $n$-sample has no measurable bias even for small sample sizes up to order 5. Thus, the only visible impact of the curvature is  on the covariance of the empirical mean.  
%
%The bias that we predict on the empirical mean in the small sample regime is difficult to verify experimentally because we have explicit formulas for geodesics only in symmetric spaces, for which. However, experiments on symmetric spaces are possible to verify our prediction on the covariance rate. 
%
To simplify further the setup and to minimize the number of parameters, we focus more particularly on constant curvature spaces, which include the Euclidean space (sectional curvature $\kappa=0$), the sphere of radius $R$ (positive sectional curvature $\kappa=1/R^2$) and the hyperbolic space of negative sectional curvature $\kappa <0$, which can be viewed as a pseudo-sphere of radius $R=-1/\sqrt{| \kappa |}$ in the Minkowski space. 
Because the space is symmetric, the expected

\subsection{Asymptotic BP-CLT for isotropic distributions}

The Hessian of the squared distance in constant curvature spaces has been computed in closed form in \cite{bhattacharya_statistics_2008} using Jacobi fields. The interested reader may also find in \cite{pennec2018} a more pedestrian approach using the embedding of the sphere (resp. the hyperbolic space) in the Euclidean space (resp. the Minkowski space) that obtains similar formulas for $\kappa = \pm 1$. 
With the notations $\log_x(y) = \theta u$ where $\theta = \| \log_x(y)\|_x$ is the distance from $x$ to $y$ and $u= \log_x(y)/\theta$ is the unit vector of $T_x\M$ pointing from $x$ to $y$, the Hessian of the squared distance is $ \frac{1}{2} H_x(y) = % f(\theta) \Id + (1-f(\theta)) u u^\top 
u u^\top  + f_{\kappa}(\theta) (\Id - u u^\top)$ with
\[
f_{\kappa}(\theta) = \left\{ \begin{array}{ll} 
\sqrt{|\kappa|}\theta \coth( \sqrt{|\kappa|} \theta ) & \text{if}\: \kappa < 0,\\
1   & \text{if}\: \kappa = 0,\\
\sqrt{\kappa}\theta \cot( \sqrt{\kappa}\theta  ) & \text{if}\: \kappa > 0.\\
\end{array} \right. .
\]
We can unify the notations for all curvatures by observing that $f_{\kappa}(\theta) = h( \kappa \theta^2 )$ where  $h(t) = \sqrt{t} \cot( \sqrt{t})$. This function is analytic at 0: its Taylor expansion is $h(t) = 1 -t/3 + O( t^2)$. Thus, the formulation is valid with positive, null and negative curvature:
\begin{equation}
\frac{1}{2} H_x(y) = % f(\theta) \Id + (1-f(\theta)) u u^\top 
u u^\top  + h( \kappa \theta^2 ) (\Id - u u^\top) 
\qquad \text{with} \qquad
h(t) = \sqrt{t} \cot( \sqrt{t})
\label{eq:HessianSquareDistanceConstantCurvature}
\end{equation}

When the point $y$ follows an isotropic  distribution at $x$ (i.e. circularly symmetric in the tangent space $T_x\M$), the distribution of $u$ is uniform on the unit sphere of dimension $d-1$, where $d$ is the dimension of the manifold. It is also independent of the distribution $dP(\theta)$ on the distance. This means that $\expect{u u^\top}= \frac{1}{d} \Id$.  
Thus, we are left with the simple expected Hessian $\bar H = 2 \gamma \Id$ with 
\begin{equation}
\gamma = \frac{1}{d} + \left(1-\frac{1}{d}\right) \overline{h}
\qquad \text{where} \qquad 
\overline{h} = \expect{h(\kappa \dist(\bar x, .)^2)} = \int_\M h ( \kappa (\dist(\bar x, y)^2 )\: \mu(dy).
\end{equation}

The BP-CLT tells us that the covariance matrix of the sample mean is:
\begin{equation}
\text{Cov}(\log_{\bar x}(\bar x_n))=  \gamma^{-2} \left( \frac{1}{n} \MM_2 \right) + O\left(\frac{1}{n^2}\right) .
\label{eq:CovarIsotropicSampleMean}
\end{equation}
Because the distribution is isotropic, the covariance matrix $\MM_2$ at the mean point $x$ is diagonal, as well as the resulting covariance on the empirical mean. Thus,  the formula boils down to a scalar equation $\text{Var}(\log_{\bar x}(\bar x_n))=  \gamma^{-2} \frac{\sigma^2}{n}  + O\left(\frac{1}{n^2}\right)$ using the variance $\sigma^2 = \text{Tr}(\MM_2)$ of the original distribution. 
In this formula, we see that the factor $\alpha = \gamma^{-2}$ plays the role of a modulation factor indicating how the rate of convergence of the covariance (or variance) differs from the Euclidean case.

\subsection{Non-asymptotic high concentration expansion}
In constant curvature spaces, the Riemannian curvature tensor is given by $R_{abcd} = \kappa\: (g_{ac}g_{db} - g_{ad}g_{cb})$. Equivalently the $(1,3)$ curvature tensor is $R^a_{bcd} = g^{ae} R_{ebcd} =   \kappa\:  (\delta^a_{c}g_{db} - \delta^a_{d}g_{cb})$.
Thus, in a normal coordinate system where $g_{ab} = \delta_{ab}$:
\[
\begin{split}
\MM_2^{cd} ( \MM_2^{ae}   R^b_{cde} + R^a_{cde} \MM_2^{be}) = &
\kappa\:  \MM_2^{cd}  \MM_2^{be} (\delta^a_{d}g_{ec} - \delta^a_{e}g_{dc}) + \kappa\:  \MM_2^{cd} (\delta^b_{d}g_{ec} - \delta^b_{e}g_{dc}) \MM_2^{ae}
\\ = & \kappa\:  ( \MM_2^{ca} g_{ec} \MM_2^{be} - \MM_2^{cd} g_{dc} \MM_2^{ba} 
+ \MM_2^{cb} g_{ec} \MM_2^{ae} - \MM_2^{cd} g_{dc} \MM_2^{ab} )
\\ = & 2\kappa ( [\MM_2 \MM_2]^{ab} - \text{Tr}(\MM_2) \MM_2^{ab} )
\end{split}
\]
Thus, Eq.(\ref{eq:EmpiricalFrechetMeanCovThm}) of Theorem \ref{eq:EmpiricalFrechetMeanBiasThm}
states that: 
\[
\text{Cov}(\bar x_n)  = \textstyle \frac{1}{n} \MM_2 \big( \Id - \frac{\kappa}{3}  \left(1-\frac{1}{n}\right) (\MM_2  - \text{Tr}(\MM_2)\Id) \big) + O\big(\epsilon^5\big)
\]
With an isotropic distribution of variance  $\sigma^2 = \text{Tr}(\MM_2)$ (i.e. $\MM_2 = {\sigma^2}/{d} \Id$), the covariance matrix of the empirical mean is diagonal and this expression boils down to:
\[
\text{Var}(\bar x_n)  = \textstyle \frac{\sigma^2}{n}   \big( 1 + \frac{2 \kappa \sigma^2}{3} \left(1-\frac{1}{d}\right)\left(1-\frac{1}{n}\right) \big)  + O\big(\epsilon^5 \big) .
\]
We can summarize the results in the following theorem.

\begin{theorem}[Modulation of the convergence speed by the curvature in space forms] 
We consider an isotropic distribution of mean $\bar x$ and variance $\sigma^2$ on a space of constant sectional curvature $\kappa$, whose support satisfies the KKC conditions. Because the space is symmetric, the bias of the empirical Fr\'echet mean $\bar x_n$ of an IID $n$-sample vanishes for small sample sizes at order 5: $\expect{\log_{\bar x}(\bar x_n)} =0 + O\big(\epsilon^5 \big) $. 

The non-asymptotic variance of the empirical Fr\'echet mean of an IID $n$-sample with sufficiently small variance $\sigma^2 < \varepsilon$ is: 
\begin{equation}
\label{eq:NonAsymptoticVarSymSpace}
\text{Var}(\bar x_n)  = \alpha \frac{\sigma^2}{n}   
\qquad\text{with} \qquad
\alpha  = \textstyle \left( 1 + \frac{ 2 }{3} \kappa \sigma^2 \left(1-\frac{1}{d}\right)\left(1-\frac{1}{n}\right) \right) + O\big(\epsilon^5 ).
\end{equation}

The asymptotic variance of the empirical Fr\'echet mean of an IID $n$-sample is:
\begin{equation}
\text{Var}(\bar x_n)  =  \alpha \frac{\sigma^2}{n}  
\qquad 
\text{with}
\qquad
\alpha = \textstyle  \left( \frac{1}{d} + \left(1-\frac{1}{d}\right) \overline{h} \right)^{-2} + O\left(\frac{1}{n^2}\right) 
,
\label{eq:AsymptoticVarSymSpace}
\end{equation}
where
\begin{equation}
\overline{h} = \expect{h(\kappa \dist(\bar x, .)^2)} \qquad \text{with} \qquad 
h(t) = \sqrt{t} \coth( \sqrt{t}) .
\end{equation}

In both expansions, the modulation factor $\alpha = \text{Var}(\bar x_n) \frac{n}{\sigma^2}$ indicates how much the variance of the empirical Fr\'echet mean deviates from the Euclidean case, i.e. how much the convergence speed is modulated by the curvature of the space: a modulation factor $\alpha >1$ indicates that the convergence is slower than in the Euclidean case, while $\alpha < 1$ indicates a faster convergence. 
\end{theorem}

A first observation is that there is no modulation for $n=1$ in the non-asymptotic high concentration expansion. Indeed, the mean of one sample is the sample itself, so that the variance of the empirical mean is the one of the sample. Note that this feature is not shared by the large sample CLT approximation. 

The second observation is that there is no modulation either in dimension 1 for both expansions, which is expected since there is no intrinsic curvature in that case. For higher dimensions, positive sectional curvature induces an increase of the dispersion of the empirical mean which slows down the convergence of the law of large number with respect to the Euclidean case. A negative sectional curvature accelerates the convergence of the law of large number with respect to the Euclidean case, as we will see more clearly below.
%In the exterme case where $\kappa \Sigma^2$ is very large (which goes out of the domain of our assumptions here), one could actually think that the speed of the convergence 
At first sight, this modulation of the rate of convergence could seem to be related to the smeary means of \cite{eltzner_smeary_2018}. However, the Taylor expansion that we get is still in $1/n$. Moreover, there is no mechanism in our expansion to create a term in $1/n^{\alpha}$ with $\alpha < 1$ in the series. The two phenomena are thus different.

\subsection{Archetypal modulation factor}
\label{sec:ArchetypalModulation}

In order to give a more intuitive idea of the impact of the curvature, we consider a uniform distribution on the Riemannian hypersphere of radius $\theta$  around the point $x$. Such a distribution maximizes the variance among isotropic distributions with support in a closed geodesic ball, and its singular distribution on $\theta$ allows us to compute the integral of the Hessian of the distance in closed form. Moreover, the covariance matrix is $\MM_2 = \theta^2/d  \Id$ so that the variance is $\theta^2$.
The modulation factor for a large number $n$ of samples is then:
\begin{equation*}
\alpha = \textstyle \left( \frac{1}{d} + \left(1-\frac{1}{d}\right) h(\kappa \theta^2) \right)^{-2} + O\left(\frac{1}{n}\right) = h(\kappa \theta^2)^{-2} + \textstyle  O(\frac{1}{d}) + O(\frac{1}{n}) .
%\label{eq:ArchetypalModulation}
\end{equation*}

\begin{wrapfigure}{r}{0.48\textwidth}
\vspace{-\baselineskip}
\includegraphics[width=0.48\textwidth]{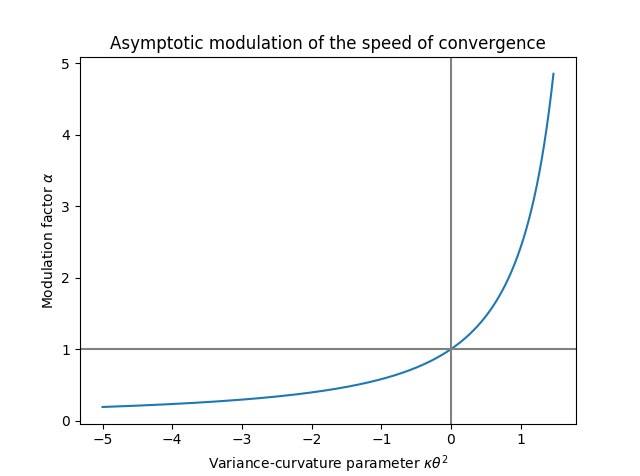}
\caption{Archetypal modulation $\tilde \alpha = h(\kappa \theta^2)^{-2}$  of the speed of convergence as a function of the curvature-variance parameter for a large dimension $d$ and a large sample $n$.   $\qquad $    $\qquad $   $\qquad  $   $\qquad $   $\qquad $   $\qquad $
\label{Fig:AsymptoticModulation}}
\end{wrapfigure}
The archetypal values of this modulation factor are obtained for positive or negative curvature with a large number of samples in a manifold of large dimension $d$: % $\tilde \alpha = h(\kappa \theta^2)^{-2}$. % + O(1/d) + O(1/n)$. 
\begin{equation}\textstyle
\tilde \alpha = h(\kappa \theta^2)^{-2} = \frac{\tan^2( \sqrt{ \kappa \theta^2} )}{ \kappa \theta^2}.
\label{eq:ArchetypalModulation}
\end{equation}
This  archetypal modulation factor is displayed in Figure (\ref{Fig:AsymptoticModulation}). The Taylor expansion of the modulation factor for a small variance shows that the modulation factor is smooth at 0 and corresponds exactly to our non-asymptotic high concentration prediction:
\[ \textstyle
\alpha= 1 + \frac{2}{3}\left( 1 -\frac{1}{d}\right) \kappa \theta^2 + \textstyle O( \theta^4, \frac{1}{n}).
\] 
It is important to notice that the variable controlling the modulation is actually $\kappa \theta^2$, i.e. the product of the sectional curvature with the variance.

For a positive variance-curvature, the modulation of the rate of convergence is above one (slower convergence than in Euclidean spaces) and actually goes to infinity when $\kappa \theta^2$ approaches $\pi^2 / 4$. Interestingly, this corresponds to the Kendall \& Karcher concentration conditions under which all this paper is restricted since we should have $\theta < \frac{\pi}{2 \sqrt{\kappa}}$ (Eq. \ref{KKCC}). Actually, as we will see below with experiments, there are examples on the sphere with $\theta = \frac{\pi}{2 \sqrt{\kappa}}$ (namely a uniform distribution on the equator) where the distribution of the empirical Fr\'echet mean converges to a mixture of Diracs rather than concentrating on a point as usual. 

For negative curvature, our formula indicates a modulation factor below 1, meaning that the variance of the estimated mean decreases faster than in Euclidean spaces. Such a phenomenon was observed in specific cases of other types of means in negatively curved manifolds \cite{basrak_limit_2010,holmes:Personnal:2015}. However, it was apparently not recognized so far as a general phenomenon of least-squares in manifolds.  For an infinite negative variance-curvature, the modulation factor actually goes to zero.  This effect suggests that we could see here the beginning of one phenomenon related to the stickiness of the Fr\'echet mean described in stratified spaces \cite{hotz_sticky_2013}. In very specific stratification cases such as corners of positive or negative curvature, one can indeed interpret the singularity with infinite curvature as the limit of a smooth manifold whose curvature is concentrating at point. In such a  process, not only the curvature is becoming (positively or negatively) infinite, but the gradient is also becoming very large, which can attract or repulse the empirical Fr\'echet mean from the singularity. 
Of course, the formula that we present here is only an approximation, so that we cannot conclude anything on this basis for large variances. However, the modulation of the convergence speed combined with the bias appear to be significant elements departing from the Euclidean situation that might partly explain the sticky mean in the range of very concentrated data on smooth affine connection or Riemannian manifolds.

\subsection{Experiments on the 2-dimensional sphere}

In order to illustrate very practically the effect of curvature on the convergence of the empirical mean estimation, we can design a very simple experiment on the sphere $\Sph_2$ embedded in $\R^3$. We consider a uniform distribution on the hypersphere of radius $\theta$ around the north pole (for $\theta \in [0,\pi]$). In the embedding 3D space, this corresponds to the horizontal small circle of radius $\sin \theta$ centered at the point $(0,0,\cos\theta)$. Such a distribution has two advantages: it is very symmetric, which considerably simplifies the computations, and it maximizes the variance among isotropic distributions with a prescribed support. Moreover, continuous or discrete distribution which entirely lie in one hemisphere (excluding the equator) are know to have a unique Fr\'echet mean located in that hemisphere \cite{buss_spherical_2001}, which is in our case at the north or the south pole. 

\begin{figure}[hbt!]
\centering
\includegraphics[width=7cm]{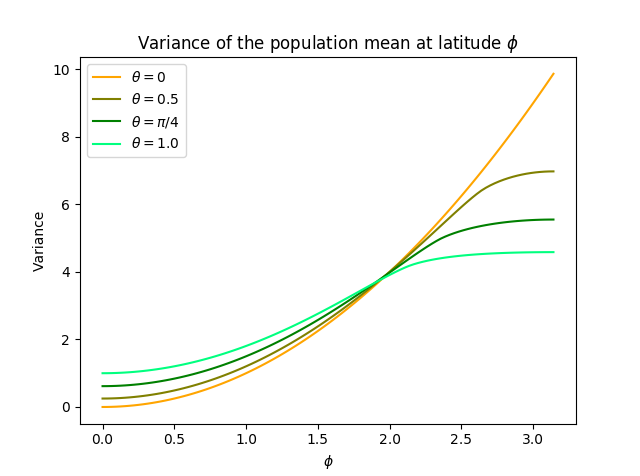}
\includegraphics[width=7cm]{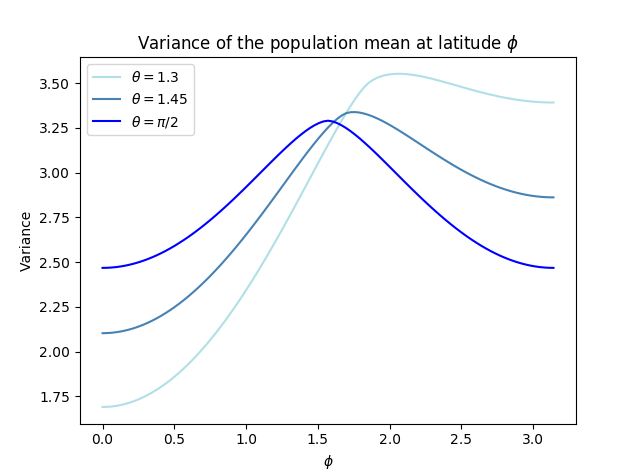}
\caption{Variance of the distribution on the circle of radius $\sin \theta$ at a point of latitude $\phi$ (counted from 0 at the north pole to $\pi$ at the south pole) on the sphere $\Sph_2$.  
\label{Fig:VarSphere}}
\end{figure}

In order to see that with simple arguments, we compute the variance at any point on the sphere.
By symmetry of our distribution, it only depends on the latitude  of the point at which we compute it. Here, we count the latitude $\phi$ from 0 at the north pole to $\pi$ at the south pole. Since the Riemannian distance is the angle between points on the sphere, the variance is:
\[
\text{Var}(\theta, \phi) =  \int_{-\pi}^{\pi} \arccos^2\left( \sin\theta \sin \phi \cos \alpha + \cos \theta \cos \phi \right) \frac{d\alpha}{2\pi} .
\]
There is generally no closed-form expression for this integral. We display in Fig.\ref{Fig:VarSphere} its value as a function of the latitude for different values of $\theta$. We clearly see that the the north pole is the unique population Fr\'echet mean for $\theta \in [0,\pi/2[$, with a variance $\text{Var}(e_3) = \theta^2$ and a covariance matrix $\text{Cov}(e_3) = \sigma^2 \Id$ with $\sigma^2 = \theta^2 / 2$.  

For low values of $\theta$, the variance is monotone and the south pole is a global maximum. 
For a value of $\theta$ larger than about 1.2 radian, a maximum appears in the south hemisphere and the south pole becomes a local minimum. For a uniform distribution on the equator ($\theta=\pi/2$),  the north and south poles become global minima with variance $\pi^2/4$, while the equator is maximizing the variance with the value $\pi^2/3$.
%\[
%\text{Var}(e_1) = \int_{-\pi}{\pi} \phi^2 d\phi/(2\pi) = 2\pi^3/3 / (2 \pi) = \pi^2/3 . 
%\]
%Thus, north and south poles are both Fr\'echet means.
For $\theta \in ]\pi/2; \pi]$, the roles of the north and south pole are  exchanged, so that we restrict in the sequel to $0 \leq \theta \leq \pi/2$.

We now draw an empirical distribution of $n$ points uniformly sampled on the horizontal small circle of radius $\sin \theta$. With the classical Gauss-Newton gradient descent algorithm $\bar x^{t+1} = \exp_{\bar x^t}( \frac{1}{n} \sum_i \log_{\bar x^t}(x_i))$, we compute the empirical Fr\'echet mean of this sample. In order to avoid the convergence problems due to potentially null gradients at maximima or saddle points, we initialize the algorithm with the north pole, which is the unique Fr\'echet mean when $\theta < \pi/2$.  Averaging the square distance of the result to the north pole for a large number $N$ of repeated random sampling allows us to compute a stochastic integral of the variance $\text{Var}(n, \theta) = \text{Var}(\bar x_n)$. Finally, we compute the normalized modulation factor $\alpha(n,\theta) = n \text{Var}(n, \theta) / \theta^2$ which indicates how the rate of convergence differs from the Euclidean case.
\begin{figure}[htb!]
\centering
\includegraphics[width=5.2cm]{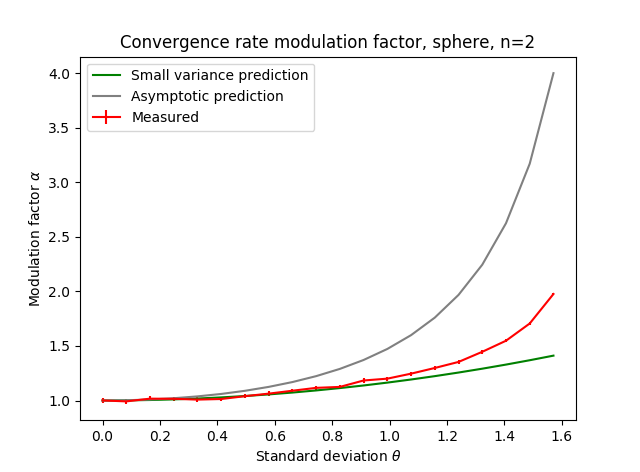}
\includegraphics[width=5.2cm]{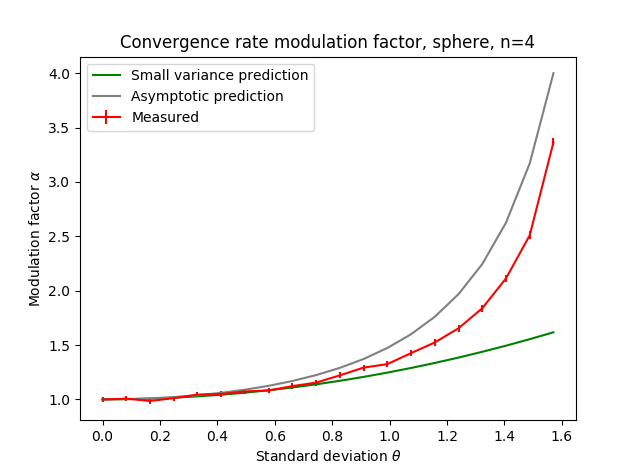}
\includegraphics[width=5.2cm]{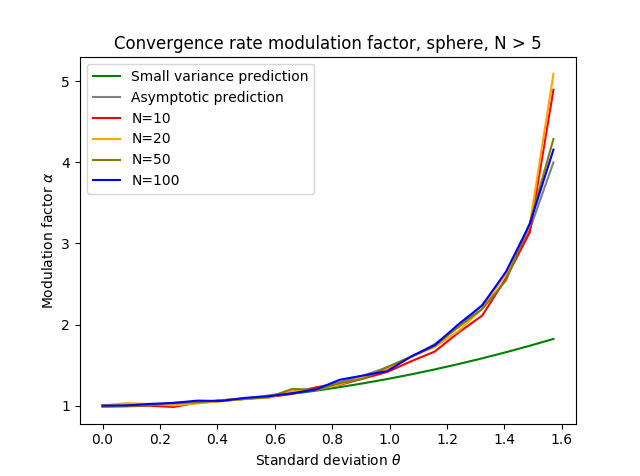}
\caption{Predicted (green curve) versus measured modulation factor $\alpha(n,\theta) = n \text{Var}(n, \theta) / \theta^2$ on the speed of convergence of the empirical Fr\'echet mean to the population Fr\'echet mean on the sphere $\Sph_2$, computed with $N=5000$ drawings to approximate the expectation for a number of samples $n=2$, $n=4$ and $n=10$ to $n=100$.  
\label{Fig:ModulationSphere}}
\end{figure}
Results are plotted in Figure \ref{Fig:ModulationSphere} for several values of $n$ and compared to the value predicted by our non-asymptotic high concentration formula Eq.(\ref{eq:NonAsymptoticVarSymSpace}) (green curve) and by the asymptotic CLT formula Eq.(\ref{eq:AsymptoticVarSymSpace}) (grey curve).
We see that the non-asymptotic high concentration expansion closely predicts the normalized modulation factor for $\theta < 0.8 rad$, whatever the number of samples. Above this value, the neglected terms in $O(\theta^5)$ take the lead and increase to a maximum which depends on how close we are to $\theta = \pi/2$. In this formula, the influence of the number $n$ of sample points is visible for a very small value, but  disappears for 5 to 10 points both in the predicted and measured modulation factors (right of Fig.\ref{Fig:ModulationSphere}).  For a small number of samples, the asymptotic CLT formula is significantly overestimating the modulation factor for a large range of values of $\theta$, while is become a very good predictor on almost all the range of values for $n \geq 10$.

\begin{figure}[htb!]
\centering
\includegraphics[width=5.2cm]{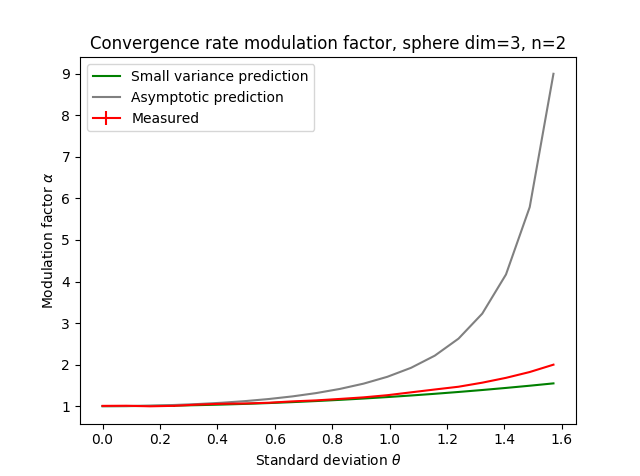}
\includegraphics[width=5.2cm]{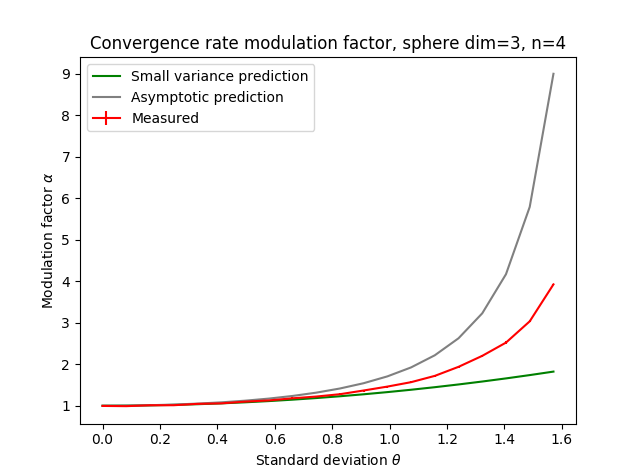}
\includegraphics[width=5.2cm]{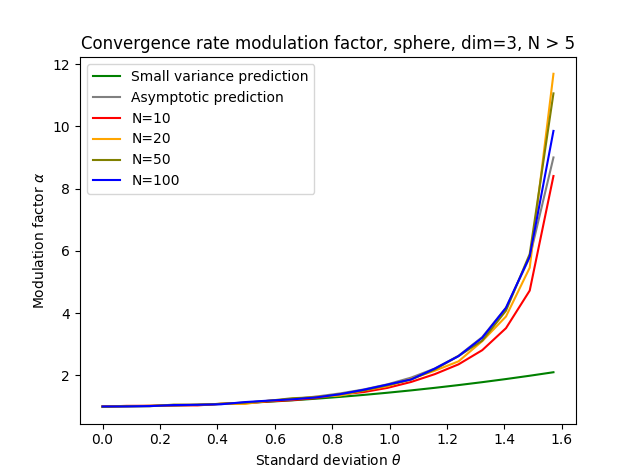}
\caption{Predicted (green curve) versus measured modulation factor $\alpha(n,\theta) = n \text{Var}(n, \theta) / \theta^2$ on the speed of convergence of the empirical Fr\'echet mean to the population Fr\'echet mean on the sphere $\Sph_3$, computed with $N=5000$ drawings to approximate the expectation for a number of samples $n=2$, $n=4$ and $n=10$ to $n=100$. 
\label{Fig:ModulationSphere3}}
\end{figure}

Very similar results are obtained for a uniform distribution on the hypersphere of radius $\theta$ on the sphere $\Sph_3$ around the north pole $(0,0,0,1)$ as shown in Fig. \ref{Fig:ModulationSphere3} and in higher dimensions. The more we increase the dimension, the more the asymptotic BP-CLT is overestimating the modulation factor for a small number of samples. However, the asymptotic approximation remains quite good for a large number of samples even close to $\theta = \pi /2$. 

\subsection{Experiments on the hyperbolic space}

A similar experiment can be made on the hyperbolic space $\Hyp_3$. We consider the positive sheet ($t>0$) of the hyperboloid of equation $x^2 + y^2 + z^2 -t^2=-1$ embedded in $\R^4$. Using the pseudo-metric $ \|(x,y,z,t)\|^2_* = x^2 + y^2 + z^2 -t^2$,  the hyperbolic space $\Hyp_3$ can be seen as the pseudo-sphere of radius -1 in the Minkowski space $\R^{3,1}$. 
With these conventions, geodesics are the trace of 2-planes passing through the origin and the Riemannian distance is the arc-length $d(x,y) = \text{arccosh}( - \scal{x}{y}_* )$. 

As in the previous case, we consider a uniform distribution on the hypersphere of radius $\theta \in \R^+$ centered at $(0,0,0,1)$ in $\Hyp_3$. In the embedding space $\R^4$, this hypersphere in $\Hyp_3$ is  the intersection  of the hyperplane $t= \cosh \theta$ with the hyperboloid. It is thus simply a sphere of radius $\sinh \theta$  in this hyperplane. A uniform distribution on the hypersphere is thus easily obtained from a uniform distribution on the 3-sphere.  Because the hyperbolic space $\Hyp_3$ has negative curvature and is Hadamard, there exists a unique Fr\'echet mean which is the center of the hypersphere $\bar x = (0,0,0,1)$. However, the negative curvature may also cause  the classical Gauss-Newton gradient descent algorithm to diverge when the variance of the distribution is large (see \cite{bini_computing_2013} for an example on SPD matrices with the affine-invariant metric). The algorithm may simply be modified with an adaptive Levenberg-Marquardt time-step $\tau < 1$ in the formula $\bar x^{t+1} = \exp_{\bar x^t}( \tau \frac{1}{n} \sum_i \log_{\bar x^t}(x_i))$.

We draw an empirical distribution of $n$ points uniformly sampled on the hypersphere of radius $\theta$.  We compute the square geodesic distance between the pole $\bar x= (0,0,0,1)$ and the obtained empirical Fr\'echet mean  $\bar x_n$of this sample.  Averaging this value for a large number $N$ of repeated random sampling allows us to compute a stochastic integral of the variance at the Fr\'echet mean $\text{Var}(n, \theta) = \text{Var}(\bar x_n)$. Finally, we compute the normalized modulation factor $\alpha(n,\theta) = n \text{Var}(n, \theta) / \theta^2$ which indicates how the rate of convergence differs from the Euclidean case.

\begin{figure}[htb!]
\centering
\includegraphics[width=5.2cm]{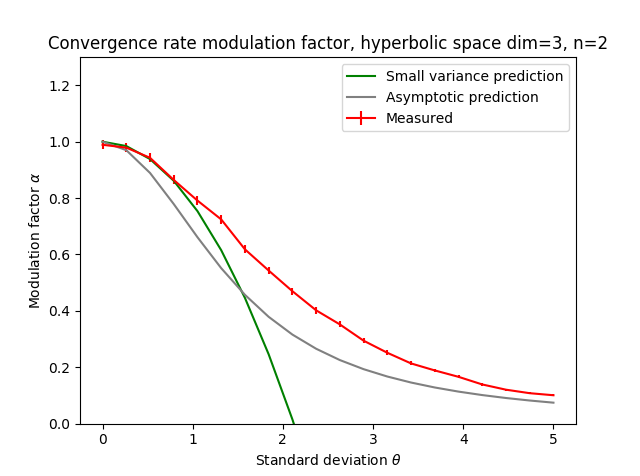}
\includegraphics[width=5.2cm]{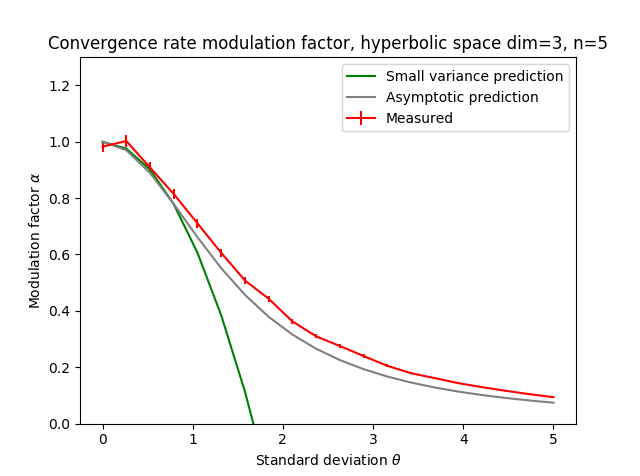}
\includegraphics[width=5.2cm]{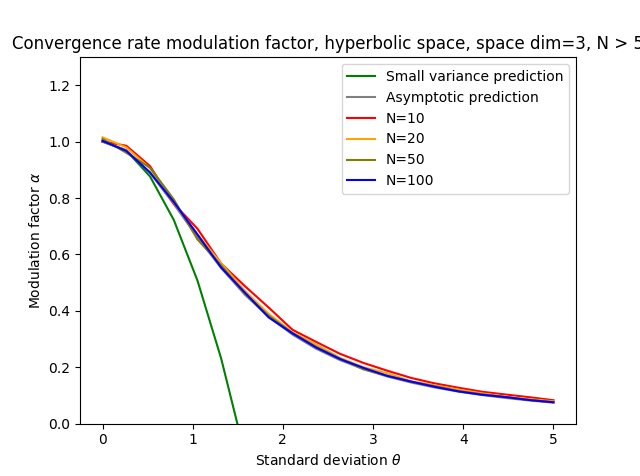}
\caption{Predicted (green curve) versus measured modulation factor $\alpha(n,\theta) = n \text{Var}(n, \theta) / \theta^2$ on the speed of convergence of the empirical Fr\'echet mean to the population Fr\'echet mean on the hyperbolic space $\Hyp_3$, computed with $N=5000$ drawings to approximate the expectation.
Top raw: curve along $\theta=[0;5]$ for a number of samples $n=2$, $n=5$ and $n=10$ to $n=100$. 
%Bottom row:  modulation factors as a fucntion of the number of samples for a fixed variance $\theta=0.1$,  
\label{Fig:ModulationHyperbolic}}
\end{figure}

Results are displayed in Figure \ref{Fig:ModulationHyperbolic} for several values of $\theta$ and $n$ (red curve) and compared to the value predicted by our non-asymptotic high concentration formula Eq.(\ref{eq:NonAsymptoticVarSymSpace}) (green curve) and by the asymptotic CLT formula Eq.(\ref{eq:AsymptoticVarSymSpace}) (grey curve).
We see that the non-asymptotic high concentration expansion closely predicts the normalized modulation factor for $\theta < 1$, whatever the number of samples. Above this value, the  terms in $O(\theta^5)$ neglected in the Taylor expansion take the lead and the formula gets useless. The asymptotic CLT formula is significantly underestimating the modulation factor for the whole range of values of $\theta$ for a small number of samples. It becomes a very good predictor on all the range of values for $n \geq 10$.
In the presented experiments, we limited the radius $\theta$ to 6 for visualization purposes.  However, we have observed a modulating factor as low as $\alpha=0.01$ for $\theta=15$ (i.e. an acceleration of the convergence by two orders of magnitude with respect to the Euclidean case) for two points. 
It is also interesting to notice that the number $n$ of points in the sample has a very limited influence on the modulation factor for $n\geq 10$: the main factor is visibly the variance $\theta^2$ of the sample.

\section{Discussion}

We have derived in this paper a new type of approximation of the moments of the empirical mean $\bar x_n$ of an IID $n$-sample of a sufficiently concentrated distribution on a Riemannian and affine connection manifolds. This high concentration expansion of the moments of the empirical mean shows a bias in $1/n$ which is directly proportional to the gradient of the curvature tensor contracted twice with the covariance matrix. This unexpected bias, apparently never described before, is important in the small sample regime. 
 The high concentration expansion of the covariance matrix of the empirical mean has also a curvature correction term modulating (accelerating or decelerating) the convergence with respect to the Euclidean case. The main variable controlling the modulation  is the covariance-curvature tensor $R(\contractvar, \contractvarcirc)\contractvar \contract \MM_2$ (the matrix $A^i_j = R^i_{klj} \MM_2^{kl}$ in coordinates). 

Our new high-concentration expansion is valid in the small sample regime. Thus, it is a useful complement to the Bhattacharya-Patrangenaru central-limit theorem, which gives an asymptotic expansion for a large number of samples. With our notations, the BP-CLT states that the empirical mean $\bar x_n$ converges in expectation to the true mean $\bar x$ with a rate of at least $1/\sqrt{n}$ and with covariance matrix (at the true mean)  $\text{Cov}(\bar x_n) = \frac{1}{n} 4 \bar H\inv \MM_2 \bar H\inv$, where $\bar H$ is the expected Hessian of the Riemannian square distance (according to the distribution).  The high concentration expansion of this expected Hessian $\bar H = 2 \Id -\frac{2}{3} R(\contractvar, \contractvarcirc)\contractvar \contract \MM_2 + O(\varepsilon^3)$ is controlled once again by the covariance-curvature tensor, as already observed in \cite{pennec2018}, so that both expansions are consistent in the high-concentration asymptotic regime. 

In constant curvature spaces, we showed that both the asymptotic BP-CLT  and the non-asymptotic high-concentration expansion predict a deviation of the decrease of the covariance of the empirical mean  with respect to the Euclidean case. This modulation of the convergence speed can be encoded with a single multiplicative factor $\alpha$ which indicates that the variance of the Fr\'echet mean decreases faster  in negatively curved space forms than in the Euclidean case, while it decreases more slowly in positive curvature space forms. The archetypal modulation factor goes to zero for an infinite negative curvature. 
This suggests that we could see here the beginning of the stickiness of the Fr\'echet mean described in stratified spaces. On the contrary, the modulation factor goes to infinity when we approach the limits of the Karcher \& Kendall concentration conditions with a uniform distribution on the equator of the sphere, for which the Fr\'echet mean is not a single point anymore. 
Although an explicit formula has previously been established for the BP-CLT in space forms, 
% \cite{bhattacharya_statistics_2008}
it seems that the interpretation of the influence of the curvature of empirical means got unnoticed so far. 

\paragraph{Acknowledgements.} This project has received funding from the European Research Council (ERC) under the European Union’s Horizon 2020 research and innovation program (grant G-Statistics agreement No 786854). Sincere thanks Yann Thanwerdas and Nicolas Guigui for their proofreading of the manuscript.

\appendix

\section{Equivalence of notations with Darling's paper}
\label{sec:Darling}

In the process of establishing the foundations of Kalman filtering on manifolds \cite{darling_geometrically_2000}, Darling derived   a coordinate free approximation of the exponential barycenter estimation in a manifold.  The work of this research report was apparently unpublished in a journal, which is very unfortunate as it got completely unnoticed.  

His Exponential Barycenter Formula Eq.(26) (Section 3.2, p.12)  actually corresponds to the order 3 of our formula  Eq.(\ref{eq:LogFrechet}). However, the 4th order necessary to see that there is a bias on the empirical mean was not derived. Rephrased in the notations of this, his formula reads 
\[
\log_x(\bar x) = \MM_1 -\frac{1}{3} R_{ijk} \MM_1^i \MM_2^{jk} + O(\varepsilon^4).
\]
However, the convention for the sign of the Riemannian curvature tensor is  the opposite of ours. Indeed, Darling defines \[
R(u,v)w = D\Gamma(v)(w\otimes u) -D\Gamma(u) (w\otimes v) + \Gamma(\Gamma(w\otimes u)\otimes v) 
- \Gamma(\Gamma(w\otimes v)\otimes u).
\] 
Taking $w=\partial_b$, $u=\partial_c$, $v=\partial_d$, we obtain 
\[
R(\partial_c, \partial_d) \partial_b = 
\partial_d \Gamma_{bc} - \partial_c\Gamma_{bd} 
+ \Gamma_{ed}\Gamma^e_{bc} - \Gamma_{ec}\Gamma^e_{bd}.
\]
Given that the connection is torsion-free, we have $\Gamma_{ij} = \Gamma_{ji}$, and we see that this formula is the opposite convention of the one used in this paper in Eq.\ref{eq:CurvatureCoord}. 

The index of the Riemannian curvature tensor also differ from ours. Here $R_{ijk} = R(\partial_i, \partial_j) \partial_k$, so that 
$R_{ijk} \MM_1^i \MM_2^{jk} = R(\MM_1,\contractvar)\contractvar \contract \MM_2$. 
Taking the opposite due to the curvature convention correctly give the term $R(\contractvar,\MM_1)\contractvar \contract \MM_2$ of our formula Eq.(\ref{eq:LogFrechet}).

\bibliographystyle{plain} 
% \bibliography{BiblioFrechetMean}

\end{document}